\documentclass[11pt]{article}
\usepackage[a4paper,left=2cm,right=2cm,top=1.5cm,bottom=2cm]{geometry}
\usepackage[utf8]{inputenc}
\usepackage[T1]{fontenc}
\usepackage[english]{babel}
\usepackage{amsthm}
\usepackage{amsmath}
\usepackage{amssymb}
\usepackage{mathrsfs}
\usepackage{setspace}
\usepackage{float}
\usepackage{multirow}
\usepackage{stmaryrd}
\usepackage{color}
\usepackage{wrapfig}
\usepackage{pifont}
\usepackage{array}
\usepackage[colorlinks=true,linkcolor=ocre,citecolor=ocre]{hyperref}
\usepackage{dsfont}
\usepackage{upgreek}
\usepackage{nicefrac}
\usepackage{tikz}{}
\usepackage{gensymb}
\usepackage{FiraSans}

\usepackage{graphicx}
\usepackage{caption}
\renewenvironment{proof}{\noindent{\sffamily{\textbf{Proof :}}}}{\begin{flushright}$\square$\end{flushright}}

\newenvironment{proofof}[2]{\noindent{\sffamily{\textbf{Proof of #1 #2 :}}}}{\begin{flushright}$\square$\end{flushright}}

\newcommand{\IE}{\mathbb{E}}
\newcommand{\IN}{\mathbb{N}}
\newcommand{\IZ}{\mathbb{Z}}

\newcommand{\IR}{\mathbb{R}}

\newcommand{\IT}{\mathbb{T}}

\newcommand{\drm}{\mathrm d}

\newcommand{\CD}{\mathcal D}

\newcommand{\CC}{\mathcal C}
\newcommand{\CH}{\mathcal H}

\newcommand{\CX}{\mathcal X}

\newcommand{\CB}{\mathcal B}

\newcommand{\CW}{\mathcal W}

\renewcommand{\P}{\mathsf{P}}
\newcommand{\PT}{\widetilde{\P}}
\newcommand{\PI}{\mathsf{\Pi}}

\newcommand{\DD}{\mathsf{D}}
\newcommand{\DC}{\mathsf{C}}

\newcommand{\DS}{\mathsf{S}}

\newcommand{\x}{\|\Xi\|_{\CX^\alpha}} 
\newcommand{\G}{H^\sharp} 

\definecolor{ocre}{RGB}{64,123,121}
\definecolor{S}{rgb}{0.0,0.5,0.0}

\newcounter{item}
\numberwithin{item}{section}

\newtheorem{theorem}[item]{\sffamily Theorem}
\newtheorem{definition}[item]{\sffamily Definition}
\newtheorem{proposition}[item]{\sffamily Proposition}
\newtheorem{lemma}[item]{\sffamily Lemma}
\newtheorem{corollary}[item]{\sffamily Corollary}

\newtheorem*{theorem*}{\sffamily Theorem}
\newtheorem*{definition*}{\sffamily Definition}
\newtheorem*{proposition*}{\sffamily Proposition}
\newtheorem*{lemma*}{\sffamily Lemma}
\newtheorem*{corollary*}{\sffamily Corollary}

\usepackage[explicit]{titlesec}
\titleformat{\section}{\centering\Large\bfseries}{\thesection \ --}{0.7em}{\Large\bfseries #1}
\titleformat{\subsection}{\centering\large\bfseries}{\thesubsection \ --}{0.4em}{\large\bfseries #1}
\titleformat{\subsubsection}{\centering\bfseries}{\thesubsubsection \ --}{0.4em}{\bfseries #1}

\let\emph\relax
\DeclareTextFontCommand{\emph}{\bfseries\em}

\providecommand{\MSC}[1]
{
	{\footnotesize	
	\textbf{MSC $\mathbf{2020}$ --} #1}
}

\providecommand{\keywords}[1]
{
	{\footnotesize	
	\textbf{Keywords --} #1}
}

\title{\bfseries Strichartz inequalities with white noise potential on compact surfaces}
\author{Antoine MOUZARD and Immanuel ZACHHUBER}
\date{}
\begin{document}

\maketitle

\abstract{We prove Strichatz inequalities for the Schrödinger equation and the wave equation with multiplicative noise on a two-dimensional manifold. This relies on the Anderson Hamiltonian described using high order paracontrolled calculus. As an application, it gives a low regularity solution theory for the associated nonlinear equations.}

\bigskip

\MSC{35J10; 60H25; 58J05}


\tableofcontents

\vspace{0.5cm}

\section*{Introduction}

Enormous progress has been made in the last decade after Hairer introduced his theory of regularity structures \cite{Hai14} and the theory of paracontrolled distributions due to Gubinelli, Imkeller, and Perkowski \cite{GIP} in the study of singular stochastic PDEs. A particular approach developed recently is the construction of a random stochastic operator to investigate associated PDEs. The first paper on this dealt with the continuum Anderson Hamiltonian, hereafter simply called Anderson Hamiltonian, and was by Allez and Chouk \cite{AllezChouk}. They used the latter theory to make sense of the operator on the two-dimensional torus, formally
$$
H=-\Delta+\xi,
$$
where $\xi$ is spatial white noise whose regularity is just below $-1$, a random field with formal covariance
$$
\IE[\xi(x)\xi(y)]=\delta_0(x-y).
$$
In the particular case of the torus, it can be constructed as the random Fourier series
$$
\xi(x)=\sum_{n\in\IZ^2}\xi_ne^{in\cdot x}
$$
with $(\xi_n)_{n\in\IZ^2}$ independent and identically distributed standard Gaussian random variables. In general, the white noise is an isometry from $L^2(M)$ to $L^2(\Omega)$ the space of random variable with finite variance. Afterwards this approach was extended to the three dimensional torus and somewhat reformulated by Gubinelli, Ugurcan and Zachhuber \cite{GUZ} and by Labbé \cite{Labbe} who used regularity structures and dealt with both periodic and Dirichlet boundary conditions. Finally, the construction was extended by Mouzard \cite{Mouzard} to the case of two-dimensional manifold using high order paracontrolled calculus.

\medskip

Naturally, substantial progress was also made in the field of singular dispersive SPDEs following the paper \cite{DebusscheWeber} due to Debussche and Weber on the cubic multiplicative stochastic Schrödinger equation and the paper \cite{GKO2} by Gubinelli, Koch and Oh on the cubic additive stochastic wave equation. Since the powerful tools from singular SPDEs are only directly applicable to parabolic and elliptic SPDEs, these initial papers were in a not so singular regime, the former using an exponential transform to remove the most singular term and the latter using a ``Da Prato-Debussche trick'' to do the same. In \cite{GUZ}, Gubinelli, Ugurcan and Zachhuber proved some sharpened results on the multiplicative Schrödinger equation and its wave analogue by reframing it in relation to the Anderson Hamiltonian as well as extending the results to dimension $3$. Moreover, Tzvetkov and Visciglia extended in \cite{TV} the results of \cite{DebusscheWeber} to a larger range of power nonlinearities, see also \cite{TVbis}. For the nonlinear wave equation with additive noise, let us mention here the follow-up paper by Gubinelli, Koch and Oh \cite{GKO3} in three dimensions with quadratic nonlinearity and the paper \cite{ORT} by Oh, Robert and Tzvetkov which extends the results of \cite{GKO2} to the case of two-dimensional surfaces and is thus salient for the current paper. 

\medskip

Let us also mention a related area of research whose aim is to solve deterministic dispersive PDEs with random initial conditions with low regularity. The study of this, which is intimately related to the analysis of invariant measures for dispersive PDEs, goes back to the seminal work of Lebowitz, Rose and Speer [25]. A series of works by Bourgain followed, let us mention here the paper [11] where a renormalization procedure similar to the current case appears but for a different reason. See also the work [14] by Burq and Tzvetkov which deals with singular random initial condition for which they obtain well-posedness results for the cubic nonlinear wave equation on a compact manifold.


\medskip

In this paper, we prove Strichartz inequalities for the Schrödinger and wave equation with white noise potential on compact surfaces. In a nutshell, Strichartz inequalities leverage dispersion in order to allow to trade integrability in time for integrability in space, see Section \ref{SectionStrichartzSchr} for a more detailed introduction and \cite{Zachhuber} where this kind of approach appeared for the Anderson Hamiltonian. Moreover, we show how this provides local well-posedness for the associated nonlinear equations in a low-regularity regimes. As for the deterministic case, the Strichartz estimates obtained depend whether the manifold has a boundary or not and are improved in the flat case of the torus. By Strichartz inequalities, we generally refer to space-time bounds on the propagators of Schrödinger and wave equations where the results on integrability are strictly better than what one gets from the Sobolev embedding so -- for definiteness we consider the Schrödinger case -- a bound like
$$
\|e^{itH}u\|_{L^p(I,L^q)}\lesssim \|u\|_{\mathcal{H}^\alpha},
$$
with $p\in[1,\infty],q>\frac{2d}{d-2\alpha}$ where $d$ denotes the dimension and $I\subset \mathbb{R}$ is an interval. The overall approach to the Schrödinger group associated to $H$ we follow is similar to the one in \cite{Zachhuber}, where such Strichartz estimates were shown for the Anderson Hamiltonian on the two and three-dimensional torus. However, one gets sharper results in the particular case of flat geometry due to the fact that one has stronger classical Strichartz inequalities available. In the more general setting of a Riemannian compact manifold, we work with a result due to Burq, Gerard and Tzvetkov \cite{BGT} which has been extended to the case with boundary by Blair, Smith and Sogge in \cite{BSSsch}. These results can be thought of as quantifying the statement ``finite frequencies travel at finite speeds -- in (frequency dependent) short time the evolution is morally on flat space''. Let us also mention at this point the recent work by Huang and Sogge \cite{HuSo} which deals with a similar setting, however their notion of singular potential refers to low integrability while in our case singular refers rather to potentials with low regularity.

\medskip

For the case of Strichartz estimates for the wave equation related to $H$, we follow the approach introduced by Burq, Lebeau and Planchon \cite{BLP} on domains with boundary. The main idea, which is why this approach is applicable, is that all that is required is that the operator driving the wave equation satisfies some growth condition on the $L^q$ bounds on the its eigenfunctions and one knows about the asymptotics of the eigenvalues, in their case the Laplace with boundary conditions. Since a Weyl law for $H$ was obtained by Mouzard in \cite{Mouzard} and our result for the Schrödinger equation gives us a suitable $L^q$ bound on the eigenfunctions of $H$, their approach turns out to be enough to prove Strichartz estimates that beat the Sobolev embedding. Overall this approach seems somewhat crude and we assume there to be sharper bounds possible whereas in the Schrödinger case, our result is the same as the one without noise obtained in \cite{BGT} worsened only by an arbitrarily small regularity loss. The state of the art of Strichartz estimates for wave equations on manifolds with boundary is the paper \cite{BSSwave}, the case of manifolds without boundary being comparable to the Strichartz estimates on Euclidean space because of the finite speed of propagation. In particular, the bounds obtained on the spectral projectors of the Anderson Hamiltonian are new and of interest in themselves.

\medskip
 
The second objective of this paper is to use the Strichartz inequalities obtained to prove local well-posedness for the associated defocussing nonlinear equations, also known as cubic multiplicative stochastic Schrödinger and wave equations. This will be done using fairly straightforward contraction arguments for which the Strichartz estimates will be crucial.

\medskip

We conclude the introduction by a brief outline of the construction of the Anderson Hamiltonian, see \cite{Mouzard} for the details. It is formally given by
$$
H=-\Delta+\xi
$$
where $\xi$ is the space white noise and belongs to $\CC^{\alpha-2}$ for any $\alpha<1$ where $\CC^\beta$ denotes the Hölder-Besov spaces recalled in Section \ref{SubsectionStrichartz}. The noise is only a distribution, rough almost surely everywhere as opposed to potential with a localised singularity, hence $Hu$ is well-defined for $u\in C^\infty$ but does not belong to $L^2$. The nature of the noise makes the naive candidate for the domain of $H$, that is the closure of
$$
\{u\in C^\infty;Hu\in L^2\}
$$
with respect to the domain norm unviable. This is precisely where the paracontrolled calculus comes into play, one can construct a random space $\CD_\Xi\subset L^2$ such that almost surely
$$
u\in\CD_\Xi\quad\implies\quad Hu\in L^2.
$$
Here $\Xi=(\xi,\Delta^{-1}\xi\cdot\xi)$ refers to the enhanced noise, see \cite{Mouzard} for its construction. The domain $\CD_\Xi$ consists of functions $u\in L^2$ paracontrolled by noise-dependent functions $X_1,X_2$ of the form
$$
u=\PT_uX_1+\PT_uX_2+u^\sharp
$$
with a remainder $u^\sharp\in\CH^2$ and the $\PT_uX_i$ are terms which are dominated by $X_i$ in terms of regularity. In particular, smooth functions do not belong to the domain in this peculiar setting. The singularity of the product is dealt with through a renormalisation procedure which corresponds to the construction of the singular term $X\cdot\xi$ where $\Delta X=\xi$. Given a regularisation $\xi_\varepsilon$ of the noise, the product $X_\varepsilon\cdot\xi_\varepsilon$ diverges but one gets a well-defined function after the subtraction of the diverging quantity
$$
c_\varepsilon:=\IE\big[X_\varepsilon\cdot\xi_\varepsilon\big].
$$
The analysis of the operator is then done with
$$
\Delta^{-1}\xi\cdot\xi:=\lim_{\varepsilon\to0}\left(\Delta^{-1}\xi_\varepsilon\cdot\xi_\varepsilon-c_\varepsilon\right),
$$
the Wick product, and the Anderson Hamiltonian corresponds to the limit of the family of operators
$$
H_\varepsilon=-\Delta+\xi_\varepsilon-c_\varepsilon
$$
as $\varepsilon$ goes to $0$. In the case of the torus, $c_\varepsilon$ is a constant due to the invariance by translation of the noise and diverges as $|\log\varepsilon|$, for details, see Section $2.1$ of \cite{Mouzard}. Note that the operator $\Delta$ is not invertible and $\Delta^{-1}$ has to be interpreted as a parametrix, that is an inverse up to a smooth term.

\medskip

While the Anderson Hamiltonian can be interpreted as the electric Laplacian $-\Delta+V$ with electric field $V=\xi$, one can consider as an analogy the magnetic Laplacian with magnetic field $B=\xi$ space white noise. This is the content of \cite{MM} where Morin and Mouzard construct
$$
H=(i\partial_1+A_1)^2+(i\partial_2+A_2)^2
$$
on the two-dimensional torus with magnetic potential $A=(A_1,A_2)$ the Lorentz gauge associated to the white noise magnetic field. Its study is also motivated by supraconductivity where $H$ plays a specific role in the third critical field of Ginzburg-Landau theory. In particular, the first results such as the self-adjointness, discrete spectrum and Weyl law holds as for the Anderson Hamiltonian, while differences are expected to appear when one looks at finer properties. Our proof for the Strichartz inequalities for the Schrödinger group associated to the Anderson Hamiltonian, which is perturbative in nature, can directly be adapted to obtain similar result for the random magnetic Laplacian with white noise magnetic field. For Strichartz estimates for the magnetic Schrödinger equation in the deterministic case, see for example \cite{DFVV} and references therein.

\bigskip

{\noindent\textbf{Organisation of the paper.}} \textit{In the first section}, we give the context for the Strichartz inequalities on manifolds in the case of the Schrödinger equations and the heat semigroup paracontrolled calculus on manifolds, see respectively Burq, Gérard, Tzvetkov \cite{BGT} and Mouzard \cite{Mouzard}. We conclude by recalling the construction of the Anderson Hamiltonian and provide new results needed in the following. \textit{In the second section}, we provide Strichartz inequalities for the Schrödinger group associated to the Anderson Hamiltonian and show how this gives local well-posedness for the stochastic cubic nonlinear Schrödinger equation with multiplicative white noise. \textit{In the last section}, we use the result on the Schrödinger group to get new bounds on the eigenvalues of the Anderson Hamiltonian and use it to prove Strichartz inequalities for the wave propagator together with the Weyl-type law. Finally, we also show how this gives local well-posedness for the stochastic cubic nonlinear wave equation with multiplicative white noise and give details for the particuliar case of the torus where one gets improved bounds.

\bigskip

\section{Preliminaries}

\subsection{Strichartz inequalities on manifolds}\label{SubsectionStrichartz}

On the torus, regularity of distributions can be measured using the Littlewood-Paley decomposition. On a manifold, one has an analogue decomposition using the eigenfunctions of the Laplace-Beltrami operator $\Delta$ as a generalisation of Fourier theory, see for example Section $2$ in \cite{ORTW} by Oh, Robert, Tzvetkov and Wang and references therein. Let $(M,g)$ be a two-dimensional compact Riemannian manifold without boundary or with boundary and Dirichlet boundary conditions. In this framework, the Laplace-Beltrami operator $-\Delta$ is a self-adjoint positive operator with discrete spectrum
$$
0\le\lambda_1<\lambda_2\le\lambda_3\le\ldots
$$ 
with the associated normalized eigenfunctions $(\varphi_n)_{n\ge1}$ belonging to $C^\infty(M)$. In the case where $M$ has no boundary, we have $\lambda_1=0$ and $\varphi_1=\textup{Vol}(M)^{-\frac{1}{2}}$ constant. Furthermore, the Weyl law gives the asymptotics
$$
\lim_{n\to\infty}\frac{\lambda_n}{n}=\frac{\textup{Vol}(M)}{4\pi}.
$$
The basis $(\varphi_n)_{n\ge1}$ of $L^2$ gives the decomposition
$$
u=\sum_{n\ge1}\langle u,\varphi_n\rangle\varphi_n
$$
for any distribution $u\in\CD'(M)$. On the torus, this gives the Littlewood-Paley decomposition of $u$ where the regularity is measured by the asymptotics behavior of $\sum_{\lambda_k\sim 2^n}\langle u,\varphi_k\rangle$. On a manifold $M$, this is done with
$$
\Delta_n:=\psi\big(-2^{-2(n+1)}\Delta\big)-\psi\big(-2^{-2n}\Delta\big)
$$
for $n\ge0$ and
$$
\Delta_{-1}:=\psi(-\Delta)
$$
with $\psi\in C_0^\infty(\IR)$ a non-negative function with $\textup{supp}(\psi)\subset[-1,1]$ and $\psi=1$ on $[-\frac{1}{2},\frac{1}{2}]$. Recall that for any function $\psi\in L^\infty(\IR)$, the operator $\psi(\Delta)$ is defined as
$$
\psi(\Delta)u=\sum_{n\ge1}\psi(\lambda_n)\langle u,\varphi_n\rangle\varphi_n
$$ 
and this yields a bounded operator from $L^2(M)$ to itself. In this setting, Besov spaces are defined for $\alpha\in\IR$ and $p,q\in[1,\infty]$ as
$$
\CB_{p,q}^\alpha:=\{u\in\CD'(M)\ ;\ \|u\|_{\CB_{p,q}^\alpha}<\infty\}
$$
where
$$
\|u\|_{\CB_{p,q}^\alpha}:=\Big(\|\Delta_{-1}u\|_{L^p(M)}^q+\sum_{n\ge0}2^{\alpha q}\|\Delta_nu\|_{L^p(M)}^q\Big)^{\frac{1}{q}}.
$$
In the particular case $p=q=\infty$ these spaces are called \textit{Hölder-Besov spaces} and we write\[
B^\alpha_{\infty,\infty}=\CC^\alpha.
\]
The case $p=q=2$ corresponds to Sobolev spaces and we have
$$
\|u\|_{\CH^\alpha}^2=\|\Delta_{-1}u\|_{L^2(M)}^2+\sum_{n\ge0}2^{2n\alpha}\|\varphi(2^{-2n}\Delta)u\|_{L^2(M)}^2
$$
where $\varphi(x):=\psi(-x^2)-\psi(-x)$.  Burq, Gérard and Tzvetkov proved in the case where $M$ has no boundary in \cite{BGT} the bound
$$
\|f\|_{L^q(M)}\lesssim\|\psi(-\Delta)f\|_{L^q(M)}+\Big(\sum_{n\ge0}\|\varphi(2^{-2n}\Delta)f\|_{L^q(M)}^2\Big)^{\frac{1}{2}}
$$
using that for $\lambda\in\IR$, we have
$$
\psi(-\lambda)+\sum_{n\ge0}\varphi(2^{-2n}\lambda)=1.
$$
Applying this to the Schrödinger group, they obtain
$$
\|e^{it\Delta}v\|_{L^p([0,1],L^q)}\lesssim\|\psi(-\Delta)v\|_{L^q(M)}+\Big\|\Big(\sum_{k\ge0}\|e^{it\Delta}\varphi(2^{-2k}\Delta)v\|_{L^q(M)}^2\Big)^{\frac{1}{2}}\Big\|_{L^p([0,1])}
$$
hence one only needs a bound for spectrally localised data. This is proved using semi-classical analysis with the use of the WKB expansion, see Proposition $2.9$ from \cite{BGT} and references therein which gives
\begin{equation}\label{Strichartzloc}
\left(\int_J\|e^{it\Delta}\varphi(h^2\Delta)v\|_{L^q(M)}^p\drm t\right)^{\frac{1}{p}}\lesssim\|v\|_{L^2(M)}
\end{equation}
for $J$ an interval of small enough length proportional to $h\in(0,1)$. Moreover, a well-known trick is to slice up the time interval into small pieces, this will be useful later. The previous bounds with the Minkowski inequality lead to
\begin{align*}
\|e^{it\Delta}v\|_{L^p([0,1],L^q)}
\lesssim\|v\|_{L^2(M)}+\Big(\sum_{k\ge0}2^{\nicefrac{2k}{p}}\|\varphi(2^{-2k}\Delta)v\|_{L^2(M)}^2\Big)^{\frac{1}{2}}\lesssim\|v\|_{\CH^{\frac{1}{p}}}.
\end{align*}
This yields the following Theorem.

\medskip

\begin{theorem}\label{thm:BGT}
Let $p\ge2$ and $q<\infty$ such that
$$
\frac{2}{p}+\frac{2}{q}=1.
$$
Then
$$
\|e^{it\Delta}u\|_{L^p([0,1],L^q)}\lesssim\|u\|_{\CH^{\frac{1}{p}}}.
$$
\end{theorem}

\medskip

While this result is optimal on general surfaces in the case $p=2$, this can be improved in the flat case of the torus. In fact, the first result concerning Strichartz inequalities for the Schrödinger equation on a compact manifold was obtain by Bourgain in \cite{Bourgain93} on the flat torus. In the case of the Anderson Hamiltonian on a compact surface without boundary we obtain the same result as Theorem \ref{thm:BGT} with an arbitrarily small loss of regularity, this is the content of Section \ref{SectionStrichartzSchr}. In the case of a surface with boundary, the following result was obtained by Blair, Smith and Sogge \cite{BSSsch}.

\medskip

\begin{theorem}\label{thm:Strichartzschboundary}
Let $M$ be a surface with boundary. Let $p\in(3,\infty]$ and $q\in[2,\infty)$ such that
$$
\frac{3}{p}+\frac{2}{q}=1
$$
Then
$$
\|e^{it\Delta}u\|_{L^p([0,1],L^q)}\lesssim\|u\|_{\CH^{\frac{2}{p}}}
$$
and 
\begin{equation}\label{Strichartzloc2}
\left(\int_J\|e^{it\Delta}\varphi(h^2\Delta)v\|_{L^q(M)}^p\drm t\right)^{\frac{1}{p}}\lesssim h^{-\frac{1}{p}}\|\varphi(h^2\Delta)v\|_{L^2(M)}
\end{equation}
for $J$ an interval of small enough length proportional to $h\in(0,1)$.
\end{theorem}

\medskip

We end this section with two classical results that will be needed in this paper. First, one still has the Bernstein Lemma with the Littlewood-Paley decomposition associated to the Laplace-Beltrami operator.

\medskip

\begin{lemma}\label{lem:bern}
Let $g:M\to \mathbb{R}$ be a function which has spectral support in an interval $[a,b]$ with $0<a<b<\infty.$ Then for any $\alpha,\beta\in \mathbb{R}$ we have the following bounds which are the analogue of Bernstein's inequality on Euclidean space
$$
\|g\|_{\mathcal{H}^\alpha}\lesssim\max (b^{\alpha-\beta},a^{\alpha-\beta}) \|g\|_{\mathcal{H}^\beta}
$$
and
$$
\|g\|_{\mathcal{H}^\alpha}\gtrsim \min (b^{\alpha-\beta},a^{\alpha-\beta}) \|g\|_{\mathcal{H}^\beta}.
$$
The former estimate still holds in the case where $a=0$ and $\alpha>\beta.$ We will chiefly apply these bounds to Littlewood-Paley projectors where $b=2a=2^j$ for $j\in\mathbb{N}$.
\end{lemma}

\medskip

\begin{proof}
The condition on $g$ means that 
$$
g=\sum_{\lambda_k\in[a,b]}(g,\phi_k)\phi_k
$$
and we have
$$
\|g\|^2_{\mathcal{H}^\alpha}=\sum_{\lambda_k\in[a,b]}(g,\phi_k)^2\lambda^{2\alpha}_k.
$$
The upper bounds follow directly with 
$$
\lambda^{2\alpha}_k=\lambda^{2\beta}_k\lambda^{2(\alpha-\beta)}_k\le \lambda^{2\beta}\max\big(b^{2(\alpha-\beta)},a^{2(\alpha-\beta)}\big)
$$ and analogously for the lower bounds.
\end{proof}

\medskip

The space $\CH^\sigma$ is an algebra only for $\sigma$ large enough depending on the dimension, this can be seen with the following Proposition and the Sobolev embedding. These types of estimates are important for the dispersive equations with cubic nonlinearity considered here.

\medskip

\begin{lemma}\label{lem:tame}
Let $\sigma\ge0$. The space $\CH^\sigma\cap L^\infty$ is an algebra and one has the bound
$$
\|f\cdot g\|_{\CH^\sigma}\lesssim\|f\|_{\CH^\sigma}\|g\|_{L^\infty}+\|g\|_{\CH^\sigma}\|f\|_{L^\infty}.
$$
\end{lemma}

\bigskip

\subsection{Basics on paracontrolled calculus}\label{SubsectionPC}

On the torus, the Littlewood-Paley decomposition can also be used to study ill-defined products. Recall that for $u\in\CD'(\IT^2)$, it is given by
$$
u=\sum_{n\ge0}\Delta_nu
$$ 
where each $\Delta_nu$ is smooth and localised in frequency in an annulus of radius $2^n$ for $n\ge1$ while the Fourier transform of $\Delta_0u$ is contained in a ball around the origin. Given two distributions $u,v\in\CD'(\IT^2)$, the product is formally given by
\begin{align*}
u\cdot v&=\sum_{n,m\ge0}\Delta_nu\cdot\Delta_mv\\
&=\sum_{n\lesssim m}\Delta_nu\cdot\Delta_mv+\sum_{n\sim m}\Delta_nu\cdot\Delta_mv+\sum_{m \lesssim n}\Delta_nu\cdot\Delta_mv\\
&=:P_uv+\Pi(u,v)+P_vu.
\end{align*}
The term $P_uv$ is called the paraproduct of $v$ by $u$ and is always well-defined while the potential singularity is encoded in the resonant term $\Pi(u,v)$. Using this decomposition, Gubinelli, Imkeller and Perkowski introduced the notion of paracontrolled calculus to develop a solution theory for singular stochastic PDEs in their seminal work \cite{GIP}; this correspond to Bony's paraproduct from \cite{Bony} in this flat case. On a manifold, an alternative paracontrolled calculus was developed by Bailleul and Bernicot in \cite{BB1} based on the heat semigroup. Instead of Littlewood-Paley, which is discrete decomposition, they used the Calderón formula
$$
u=\lim_{t\to0}P_tu=\int_0^1Q_tu\frac{\drm t}{t}+P_1u
$$ 
for $u\in\CD'(M)$ with $P_t$ the heat semigroup and $Q_t=-t\partial_tP_t$. Using Gaussian upper bounds for the heat kernel and its derivatives, this defines a continuous analogue of the Littlewood-Paley decomposition where $\sqrt{t}\simeq 2^{-n}$ and yields descriptions of Besov-Hölder and Sobolev spaces for scalar fields on manifolds. This can be used to construct a paraproduct $\P$ and a resonant product $\PI$ such that
$$
u\cdot v=\P_uv+\PI(u,v)+\P_vu
$$
that verify the same important properties of their Fourier analogue $P$ and $\Pi$. It was later extended to a higher order paracontrolled calculus by Bailleul, Bernicot and Frey in \cite{BB2} to deal with rougher noise than the initial work of Gubinelli, Imkeller and Perkowski, again in a general geometric framework. While these different works dealt with parabolic PDEs, the paracontrolled calculus can be used to study singular random operators. It was first used by Allez and Chouk in \cite{AllezChouk} to study the Anderson Hamiltonian
$$
H=-\Delta+\xi
$$
on the two-dimensional torus. The same operator was constructed on the torus by Gubinelli, Ugurcan and Zachhuber in \cite{GUZ} on $\IT^d$ with $d\in\{2,3\}$ to solve associated evolution PDEs. Labbé also constructed in \cite{Labbe} the operator in two and three dimensions with different boundary conditions using regularity structures. Finally, Mouzard used the heat semigroup paracontrolled calculus in \cite{Mouzard} to construct the operator on a two-dimensional manifold and obtained an almost sure Weyl-type law. Note that the work \cite{Mouzard} is self-contained and a gentle introduction to the paracontrolled calculus on manifolds in the spatial framework. For another example of singular random operators, see \cite{MM} where Morin and Mouzard construct the magnetic Laplacian with white noise magnetic field on $\IT^2$. With this work, we show that this approach is also well-suited for the study of dispersive PDEs.

\medskip

The heat semigroup paracontrolled calculus is a theory to study PDEs with singular products on manifolds. Given a suitable family $(V_i)_{1\le i\le d}$ of first order differential operators, one can construct a paraproduct $\P$ and a resonant term $\PI$ based on the heat semigroup associated to
$$
L:=-\sum_{i=1}^dV_i^2.
$$
We briefly outline this construction here, see \cite{BB2,BB3} in the parabolic space-time setting and \cite{Mouzard} in the space setting for the details. In particular, the Laplace-Beltrami operator on a manifold can be written in this form, see for example Stroock’s book \cite{Stroock}. For any distribution $u\in\CD'(M)$, the heat semigroup
$$
P_tu=e^{-tL}u
$$
provides a smooth approximation as $t$ goes to $0$. Introducing its derivative
$$
Q_t:=-t\partial_tP_t,
$$
one gets an analogue of the Littlewood-Paley decomposition as explained before. While the $\Delta_n$'s enjoy proper orthogonal relation in the sense that $\Delta_n\Delta_m$ is equal to zero for $|n-m|>1$, we only have in this continuous framework
$$
Q_tQ_s=\frac{ts}{(t+s)^2}\big((t+s)L\big)^2e^{-tL}
$$
which is indeed small if $s\ll t$ or $t\ll s$. For a given integer $b\in\IN^*$, let
$$
Q_t^{(b)}:=(tL)^be^{-tL}.
$$
Then
$$
Q_t^{(b)}Q_s^{(b)}=\left(\frac{ts}{(t+s)^2}\right)^bQ_{t+s}^{(b)}
$$
hence the parameter $b$ encodes a cancellation properties between different scales $t$ and $s$. Furthermore, we have
$$
\int_0^1Q_t^{(b)}u\frac{\drm t}{t}=\lim_{t\to 0}P_t^{(b)}u=u
$$
where $P_0^{(b)}=\textup{Id}$ and
$$
-t\partial_tP_t^{(b)}=Q_t^{(b)}.
$$
In  particular, we have $P_t^{(b)}=p_b(tL)e^{-tL}$ with $p_b$ a polynomial of degree $b-1$ such that $p_b(0)=1$. Denote by $\mathsf{StGC}^a$ the family of operators $(Q_t)_{t\in[0,1]}$ of the form
$$
Q_t=(t^{\frac{|I|}{2}}V_I)(tL)^je^{-tL}
$$
with $a=|I|+2j$ and $\mathsf{GC}^a$ the operator with kernel satisfying Gaussian upper bounds with cancellation of order $a$, see Section $1.2$ \cite{Mouzard} for the definitions. We have
\begin{align*}
u\cdot v&=\lim_{t\to 0}P_t^{(b)}\big(P_t^{(b)}u\cdot P_t^{(b)}v\big)\\
&=\int_0^1Q_t^{(b)}\big(P_t^{(b)}u\cdot P_t^{(b)}v\big)\frac{\drm t}{t}+\int_0^1P_t^{(b)}\big(Q_t^{(b)}u\cdot P_t^{(b)}v\big)\frac{\drm t}{t}+\int_0^1P_t^{(b)}\big(P_t^{(b)}u\cdot Q_t^{(b)}v\big)\frac{\drm t}{t}\\
&\quad+P_1^{(b)}\big(P_1^{(b)}u\cdot P_1^{(b)}v\big).
\end{align*}
After a number of integrations by parts, we get
$$
u\cdot v=\P_uv+\PI(u,v)+\P_vu
$$
where $\P_uv$ is a linear combination of terms of the form
$$
\int_0^1Q_t^{1\bullet}\big(P_tu\cdot Q_t^2v\big)\frac{\drm t}{t}
$$
and $\PI(u,v)$ of
$$
\int_0^1P_t^{\bullet}\big(Q_t^1u\cdot Q_t^2v\big)\frac{\drm t}{t}
$$
where $Q^1,Q^2\in\mathsf{StGC}^{\frac{b}{2}}$ and $P\in\mathsf{StGC}^{[0,b]}$. In general, the operator $V_i$'s do not commute hence the need for the notation
$$
Q_t^\bullet=\Big((t^{\frac{|I|}{2}}V_I)(tL)^je^{-tL}\Big)^\bullet:=(tL)^je^{-tL}(t^{\frac{|I|}{2}}V_I)
$$
which comes from the integration by parts.  For simplicity we state most of the results of this Section in Besov-Hölder spaces. The following Proposition gives the continuity estimates of the paraproduct and the resonant term between Sobolev and Hölder-Besov functions but they hold in the same way by replacing all the Sobolev spaces by Besov-Hölder spaces.

\medskip

\begin{proposition}\label{Pbounds}
Let $\alpha,\beta\in(-2b,2b)$ be regularity exponents.
\begin{itemize}
	\item[$\centerdot$] If $\alpha\ge0$, then $(f,g)\mapsto\P_fg$ is continuous from $\mathcal{H}^\alpha\times\CC^\beta$ to $\CH^\beta$.
	\item[$\centerdot$] If $\alpha<0$, then $(f,g)\mapsto\P_fg$ is continuous from $\CH^\alpha\times\CC^\beta$ to $\CH^{\alpha+\beta}$.
	\item[$\centerdot$] If $\alpha+\beta>0$, then $(f,g)\mapsto\PI(f,g)$ is continuous from $\CH^\alpha\times\CC^\beta$ to $\CH^{\alpha+\beta}$.
\end{itemize}
\end{proposition}

\medskip

While $\P$ and $\PI$ are tools to describe products, the interwined paraproduct $\PT$ naturally appears when formulating solutions to PDEs. The intertwining relation is
$$
L\circ\PT=\P\circ L
$$
hence $\PT_uv$ is given as a linear combination of
\begin{align*}
\int_0^1L^{-1}Q_t^{1\bullet}\big(P_tu\cdot Q_t^2Lv\big)\frac{\drm t}{t}&\sim\int_0^1(tL)^{-1}Q_t^{1\bullet}\big(P_tu\cdot Q_t^2(tL)v\big)\frac{\drm t}{t}\\
&\sim\int_0^1\widetilde Q_t^{1\bullet}\big(P_tu\cdot \widetilde Q_t^2v\big)\frac{\drm t}{t}
\end{align*}
with $\widetilde Q^1\in\mathsf{StGC}^{\frac{b}{2}-2}$ and $\widetilde Q^2\in\mathsf{StGC}^{\frac{b}{2}+2}$. The operator $L$ is not invertible and everything here is done up to a smooth error term, see \cite{Mouzard}. In particular, $\PT$ has the same structure as $\P$ for large $b$ and satisfies the same continuity estimates as $\P$. Intuitively, the intertwined operator $\PT$ describes solutions to elliptic PDEs of the form
$$
Lu=u\xi=\P_u\xi+\P_\xi u+\PI(u,\xi)
$$
which rewrites
$$
u=\PT_u(L^{-1}\xi)+u^\sharp
$$
hence this is the operator used to described the domain $\CD_\Xi$ of the Anderson Hamiltonian. The final ingredient of paracontrolled calculus is a toolbox of correctors and commutators made to express the singular product between a paracontrolled functions $u$ and the noise $\xi$ in a form involving only ill-defined expressions of the noise independent of $u$. The first one introduced by \cite{GIP} is in this framework the corrector
$$
\DC(u,X,\xi):=\PI\big(\PT_uX,\xi\big)-u\PI(X,\xi)
$$
which translates the rough paths philosophy: the multiplication of a function that locally looks like $X$ with $\xi$ is possible if one is given the multiplication of $X$ itself with $\xi$. This is the content of the following Proposition.

\medskip

\begin{proposition}\label{Cbound}
Let $\alpha\in(0,1)$ and $\beta,\gamma\in\IR$. If
$$
\alpha+\beta<0\quad\text{and}\quad\alpha+\beta+\gamma>0,
$$
then $\DC$ extends in a unique continuous operator from $\CC^\alpha\times\CC^\beta\times\CC^\gamma$ to $\CC^{\alpha+\beta+\gamma}$.
\end{proposition}

\medskip

While we do not give the proof, one has the following heuristic. For any $x\in M$, we have
\begin{align*}
\DC\big(f,g,h\big)(x)&=\PI\big(\PT_fg,h\big)(x)-f(x)\cdot\PI\big(g,f\big)(x)\\
&=\PI\big(\PT_fg-f(x)\cdot g,h\big)(x)\\
&\simeq\PI\big(\PT_{f-f(x)}g,h\big)(x)
\end{align*}
where $\simeq$ is equal up to a smooth term since $g\simeq\PT_1g$. Since $f\in\CC^\alpha$ with $\alpha\in(0,1)$, the term $f-f(x)$ allows to gain regularity in the paraproduct using that $\alpha+\beta<0$ ending up with a term of better regularity $\alpha+\beta+\gamma>0$. Continuity results on a number of correctors and commutators and their iterated version are also available, we refer to \cite{Mouzard} and references therein for further details. For example, one needs the swap operator
$$
\DS(f,g,h)=\P_h\PT_fg-\P_f\P_hg
$$
for the study of the Anderson Hamiltonian which is continuous from $\CH^\alpha\times\CC^\beta\times\CC^\gamma$ to $\CH^{\alpha+\beta+\gamma}$ for $\alpha,\beta\in\IR$ and $\gamma<0$.

\bigskip

\subsection{Construction of the Anderson Hamiltonian}\label{SubsectionAnderson}

In this Section, we recall the ideas behind the construction of the Anderson Hamiltonian with the heat semigroup paracontrolled calculus as done in \cite{Mouzard} and state the important results we shall use without proofs. We also provide new straightforward results from the construction needed for our proof of Strichartz inequalities. The Anderson Hamiltonian on a two-dimensional manifold $M$ is formally given by
$$
H:=L+\xi
$$
where $-L$ is the Laplace-Beltrami operator and $\xi$ is a spatial white noise. The noise belongs almost surely to $\CC^{\alpha-2}$ for any $\alpha<1$ hence the product of $\xi$ with a generic $L^2$-function in not defined almost surely. As explained, it was first constructed by Allez and Chouk in \cite{AllezChouk} on $\IT^2$. We work here with the construction on a two-dimensional manifolds from \cite{Mouzard} using the high order paracontrolled calculus since this is the setting in which we want to prove Strichartz inequalities. Following the recent development in singular stochastic PDEs, the idea is to construct a random almost surely dense subspace $\CD_\Xi$ of $L^2$ such that the operator makes sense for $u\in\CD_\Xi\subset L^2$ with $\Xi$ an enhancement of the noise that depend only measurably on the noise $\xi$. One can then prove that $H$ is self-adjoint with discrete spectrum
$$
\lambda_1(\Xi)\le\lambda_2(\Xi)\le\ldots\le\lambda_n(\Xi)\le\ldots 
$$
and compare it to the eigenvalues of the Laplace-Beltrami operator $(\lambda_n)_{n\ge1}$. While the construction of the domain $\CD_\Xi$ relied on the notion of strongly paracontrolled functions in \cite{AllezChouk,GUZ}, the high order paracontrolled calculus gives a finer description of the domain. In particular, it yields sharp bounds on the eigenvalues of the form
$$
\lambda_n-m_\delta^1(\Xi)\le\lambda_n(\Xi)\le(1+\delta)\lambda_n+m_\delta^2(\Xi)
$$
for any $\delta\in(0,1)$ and $m_\delta^1(\Xi),m_\delta^2(\Xi)>0$ random constants depending on the enhanced noise $\Xi$, see \cite{Mouzard} for a precise construction. In particular, it implies the almost sure Weyl-type law
$$
\lim_{\lambda\to\infty}\lambda^{-1}\big|\{n\ge0;\lambda_n(\Xi)\le\lambda\}\big|=\frac{\textup{Vol}(M)}{4\pi}.
$$  
We briefly present the construction of $H$ and refer to \cite{Mouzard} for the details.

\medskip

Coming from Lyons' rough paths \cite{Lyons} and Gubinelli's controlled paths \cite{Gub} which were developed as a pathwise approach to stochastic integration, the method used over the last decade to solve singular stochastic PDEs is to work in random subspaces of classical function spaces built from the noise tailor-made for the problem under consideration. In the context of singular random operators, this corresponds to the construction of a random dense domain $\CD_\Xi\subset L^2$ on which the operator almost surely makes sense. In the framework of paracontrolled calculus, one considers functions $u$ paracontrolled by noise-dependent reference functions of the form
$$
u=\PT_{u'}X+u^\sharp
$$
where the new unknown is $(u',u^\sharp)$. The function $u'$ has to be thought as the ``derivative'' of $u$ with respect to $X$ while the error $u^\sharp$ is a smoother remainder. The goal is to find a paracontrolled expression for $u\in L^2$ such that $Hu\in L^2$. Let us first assume that $u$ is smooth, then we formally get
\begin{align*}
Lu&=Hu-u\xi\\
&=-\P_u\xi+Hu-\P_\xi u-\PI(u,\xi)\in\CH^{\alpha-2+\kappa}
\end{align*}
for any $\kappa>0$. Indeed, the term of lowest regularity is the paraproduct $\P_u\xi\in\CH^{\alpha-2}$ since $u\in L^2$ and $\xi\in\CC^{\alpha-2}$. Then elliptic regularity theory gives $u\in\CH^\alpha$ and suggests for $u$ the paracontrolled form
$$
u=\PT_uX+u^\sharp
$$
with $X=-L^{-1}\xi$ and $u^\sharp\in\CH^{2\alpha}$. Given such a function, the resonance between $u$ and $\xi$ can be described by
\begin{align*}
\PI(u,\xi)&=\PI\big(\PT_uX,\xi\big)+\PI(u^\sharp,\xi)\\
&=u\PI(X,\xi)+\DC(u,X,\xi)+\PI(u^\sharp,\xi)
\end{align*}
using the corrector $\DC$ since $\PI\big(\PT_uX,\xi\big)$ is not defined due to lack of regularity, see Propositions \ref{Pbounds} and \ref{Cbound}. Since $3\alpha+2>0$, the only term on the right hand side which is potentially undefined is $\PI(X,\xi)$ and its definition is independent of the study of $H$, see \cite{Mouzard} for more details including the renormalisation. Given the enhanced data 
$$
\Xi:=\big(\xi,\PI(X,\xi)\big)\in\CC^{\alpha-2}\times\CC^{2\alpha-2}=:\CX^\alpha
$$ 
one can define the Anderson Hamiltonian $H$ on
$$
\CD:=\{u\in L^2;u-\PT_uX\in\CH^{2\alpha}\}\subset\CH^\alpha
$$
with
$$
Hu:=Lu+\P_\xi u+u\PI(X,\xi)+\DC(u,X,\xi)+\PI(u^\sharp,\xi).
$$
However, this gives only an unbounded operator $(H,\CD)$ from $\CH^\alpha\subset L^2$ to $\CH^{2\alpha-2}$ which is not a subspace of $L^2$ and thus $H$ will not takes value in $L^2$ a priori. A finer description of the domain with a second order paracontrolled expansion allows to construct a dense subspace $\CD_\Xi\subset L^2$ such that $(H,\CD_\Xi)$ is an unbounded operator on $L^2$. Using the classical theory for unbounded operators, it is possible to prove that $H$ is self-adjoint with pure point spectrum. In the expression for $H$, the roughest term is
$$
\P_\xi u+\P_u\PI(X,\xi)\in\CH^{2\alpha-2}.
$$
To cancel it with a paracontrolled expansion, we use the commutator $\DS$ to get
\begin{align*}
\P_\xi u&=\P_\xi\PT_uX+\P_\xi u^\sharp\\
&=\P_u\P_\xi X+\DS(u,X,\xi)+\P_\xi u^\sharp
\end{align*}
hence the roughest term is
$$
\P_u\P_\xi X+\P_u\PI(X,\xi)\in\CH^{2\alpha-2}.
$$
In the end, it is cancelled with the paracontrolled expansion
$$
u=\PT_uX_1+\PT_uX_2+u^\sharp
$$
where
$$
X_1:=-L^{-1}\xi\quad\text{and}\quad X_2:=-L^{-1}\big(\P_\xi X_1+\PI(X_1,\xi)\big).
$$

\medskip

\begin{definition}
We define the space $\CD_\Xi$ of functions paracontrolled by $\Xi$ as
$$
\CD_\Xi:=\big\{u\in L^2;\ u^\sharp:=u-\PT_uX_1-\PT_uX_2\in\CH^2\big\}.
$$
\end{definition}

\medskip

A powerful tool to investigate the domain $\CD_\Xi$ and $H$ is the $\Gamma$ map defined as follows. The domain is given as
$$
\CD_\Xi=\Phi^{-1}(\CH^2)
$$
with
$$
\Phi(u):=u-\PT_u(X_1+X_2).
$$
The map $\Phi$ is not necessarily invertible so we introduce a parameter $s>0$ and consider the map
$$
\Phi^s(u):=u-\PT_u^s(X_1+X_2)
$$
where $\PT^s$ is a truncated paraproduct. In particular, $\PT^s$ goes to $0$ as $s$ goes to $0$ and the difference $\PT-\PT^s$ is smooth for any $s>0$. This has to be thought as a frequency cut-off where one gets rid of a number of low frequencies in order to make a term small. Thus $\Phi^s$ is a perturbation of the identity of $\CH^\beta$ for any $\beta\in[0,\alpha)$ and thus invertible for $s=s(\Xi)$ small enough. We define $\Gamma$ to be its inverse which is implicitly defined by
$$
\Gamma u^\sharp=\PT_{\Gamma u^\sharp}^s(X_1+X_2)+u^\sharp
$$
for any $u^\sharp\in\CH^\beta$. It will be a crucial tool to describe the operator $H$ since
$$
\CD_\Xi=\Phi^{-1}(\CH^2)=(\Phi^s)^{-1}(\CH^2)=\Gamma(\CH^2)
$$
where the equality holds because the difference $\PT-\PT^s$ is smooth. Of course the map $\Gamma$ depends on the choice of $s$, however the above reasoning tells us that the image of $\Gamma$ does not change by changing $s$ so we omit this dependence in the sequel. The maps $\Phi^s$ and $\Gamma$ satisfy a number of continuity estimates that we shall use throughout this work, this is the content of the following Proposition. Let
$$
s_\beta(\Xi):=\left(\frac{\alpha-\beta}{m\x(1+\x)}\right)^{\frac{4}{\alpha-\beta}}
$$
for any $0\le\beta<\alpha$. Note that the bounds in Sobolev and Hölder spaces are proved directly while the bounds in $L^p$ follow by interpolation as in \cite{Zachhuber}.

\medskip

\begin{proposition}\label{prop:gamma}
Let $\beta\in[0,\alpha)$ and $s\in(0,1)$. We have
$$
\|\Phi^s(u)-u\|_{\CH^\beta}\le\frac{m}{\alpha-\beta}s^{\frac{\alpha-\beta}{4}}\x(1+\x)\|u\|_{L^2}.
$$
If moreover $s<s_\beta(\Xi)$, this implies
$$
\|\Gamma u^\sharp\|_{\CH^\beta}\le\frac{1}{1-\frac{m}{\alpha-\beta}s^{\frac{\alpha-\beta}{4}}\x(1+\x)}\|u^\sharp\|_{\CH^\beta}
$$
as well as the same bounds in $\CC^\beta$. The map $\Phi$ is  also continuous from $L^p$ to itself for $p\in[1,\infty]$ and $\CH^\sigma$ to itself for $\sigma\in[0,1)$ while the same holds for $\Gamma$ provided $s$ is small enough.
\end{proposition}

\medskip

Let us insist that the norm $\CH^\beta$ of $u_s^\sharp:=\Phi^s(u)$ is always controlled by $\|u\|_{\CH^\beta}$ while $s$ needs to be small depending on the noise for $\|u\|_{\CH^\beta}$ to be controlled by $\|u_s^\sharp\|_{\CH^\beta}$. We also define the map $\Gamma_\varepsilon$ associated to the regularised noise $\Xi_\varepsilon$ as
$$
\Gamma_\varepsilon u^\sharp=\PT_{\Gamma_\varepsilon u^\sharp}^sX_1^{(\varepsilon)}+\PT_{\Gamma_\varepsilon u^\sharp}^sX_2^{(\varepsilon)}+u^\sharp
$$
with 
$$
-LX_1^{(\varepsilon)}:=\xi_\varepsilon\quad\text{and}\quad-LX_2^{(\varepsilon)}:=\PI(X_1^{(\varepsilon)},\xi_\varepsilon)-c_\varepsilon+\P_{\xi_\varepsilon}X_1^{(\varepsilon)}.
$$ 
It satisfies the same bounds as $\Gamma$ with constants which depend in an increasing way on $\|\Xi_\varepsilon\|_{\CX^\alpha}\lesssim 1+\|\Xi\|_{\CX^\alpha}$ and the following approximation Lemma holds. Thus we may choose $s$ independently of $\varepsilon$.

\medskip

\begin{lemma}\label{GammaConvergence}
For any $0\le\beta<\alpha$ and $0<s<s_\beta(\Xi)$, we have
$$
\|\textup{Id}-\Gamma\Gamma_\varepsilon^{-1}\|_{L^2\to\CH^\beta}\lesssim_{\Xi,s,\beta}\|\Xi-\Xi_\varepsilon\|_{\CX^\alpha}.
$$
In particular, this implies the norm convergence of $\Gamma_\varepsilon$ to $\Gamma$ with the bound
$$
\|\Gamma-\Gamma_\varepsilon\|_{\CH^\beta\to\CH^\beta}\lesssim_{\Xi,s,\beta}\|\Xi-\Xi_\varepsilon\|_{\CX^\alpha}.
$$
\end{lemma}

\medskip

In particular, this allows to prove density of the domain.

\medskip

\begin{corollary*}
The domain $\CD_\Xi$ is dense in $\CH^\beta$ for any $\beta\in[0,\alpha)$.
\end{corollary*}

\medskip

For any $u\in\CD_\Xi$, the operator $H$ is given by
$$
Hu=Lu^\sharp+\P_\xi u^\sharp+\PI(u^\sharp,\xi)+R(u)
$$
with $u^\sharp=\Phi(u)\in\CH^2$ and $R$ an explicit operator depending on $\Xi$ which is continuous from $\CH^\alpha$ to $\CH^{3\alpha-2}$. For each $s>0$, we have a different representation of $H,$ namely
$$
Hu=H\Gamma u_s^\sharp=Lu_s^\sharp+\P_\xi u_s^\sharp+\PI(u_s^\sharp,\xi)+R(\Gamma u_s^\sharp)+\Psi^s(\Gamma u_s^\sharp)
$$
with $u_s^\sharp=\Phi^s(u)\in\CH^2$ and $\Psi^s$ an explicit operator depending on $\Xi$ and $s$ continuous from $L^2$ to $C^\infty$ which we henceforth include in the operator $R$. The operator $H\Gamma$ is thus a perturbation of $L$, the following Proposition shows that it is a continuous operator from $\CH^2$ to $L^2$. In Section \ref{SectionStrichartzSchr}, we show that it is even a lower order perturbation of the Laplace-Beltrami operator; this will be crucial to obtain Strichartz inequalities.

\medskip

\begin{proposition}\label{SobolevBoundH}
For any $\gamma\in(-\alpha,3\alpha-2)$ and $s$ as above, we have
$$
\|Hu\|_{\CH^\gamma}=\|H\Gamma u_s^\sharp\|_{\CH^\gamma}\lesssim\|u_s^\sharp\|_{\CH^{\gamma+2}}
$$
with $u=\Gamma u_s^\sharp\in\CD_\Xi$. In particular, the result holds for $\gamma\in(-1,1)$ since the noise belongs to $\CC^{\alpha-2}$ for any $\alpha<1$.
\end{proposition}

\medskip

\begin{proof}
We have
$$
H\Gamma u_s^\sharp=Lu_s^\sharp+\P_\xi u_s^\sharp+\PI(u_s^\sharp,\xi)+R(u)
$$
with $u=\Gamma u_s^\sharp$. Assume first that $0<\gamma<3\alpha-2$ hence
\begin{align*}
\|H\Gamma u_s^\sharp\|_{\CH^\gamma}&\lesssim\|Lu_s^\sharp\|_{\CH^\gamma}+\|\P_\xi u_s^\sharp+\PI(u_s^\sharp,\xi)\|_{\CH^\gamma}+\|R(u)\|_{\CH^\gamma}\\
&\lesssim\|u_s^\sharp\|_{\CH^{\gamma+2}}+\|\xi\|_{\CC^{\alpha-2}}\|u_s^\sharp\|_{\CH^{\gamma+2-\alpha}}+\|R(u)\|_{\CH^{3\alpha-2}}
\end{align*}
where the condition $\gamma>0$ is needed for the resonant term and $\gamma<3\alpha-2$ for $R(u)$. The result follows for this case since
$$
\|R(u)\|_{\CH^{3\alpha-2}}\lesssim\|u\|_{\CH^\alpha}\lesssim\|u_s^\sharp\|_{\CH^\alpha}\lesssim\|u_s^\sharp\|_{\CH^{\gamma+2}}.
$$
Assume now that $-\alpha<\gamma\le0$. For any $\delta>0$, we have
\begin{align*}
\|H\Gamma u_s^\sharp\|_{\CH^\gamma}&\lesssim\|Lu_s^\sharp\|_{\CH^\gamma}+\|\P_\xi u_s^\sharp+\PI(u_s^\sharp,\xi)\|_{\CH^\gamma}+\|R(u)\|_{\CH^\gamma}\\
&\lesssim\|Lu_s^\sharp\|_{\CH^\gamma}+\|\P_\xi u_s^\sharp+\PI(u_s^\sharp,\xi)\|_{\CH^\delta}+\|R(u)\|_{\CH^\gamma}\\
&\lesssim\|u_s^\sharp\|_{\CH^{\gamma+2}}+\|\xi\|_{\CC^{\alpha-2}}\|u_s^\sharp\|_{\CH^{\delta+2-\alpha}}+\|R(u)\|_{\CH^{3\alpha-2}}
\end{align*}
using that $\gamma\le0<\delta$. The proof is complete since $\gamma>-\alpha$ and $\delta$ small enough implies $\gamma+2>\delta+2-\alpha$.
\end{proof}

\medskip

As the parameter $s>0$ yields different representation of $H$, the domain $\CD_\Xi$ is naturally equipped with the norms
$$
\|u\|_{\CD_\Xi}^2:=\|u\|_{L^2}^2+\|u_s^\sharp\|_{\CH^2}^2.
$$
which are equivalent to the graph norm
$$
\|u\|_{H}^2:=\|u\|_{L^2}^2+\|Hu\|_{L^2}^2.
$$
In particular, this shows that the operator $H$ is closed on its domain $\CD_\Xi$.

\medskip

\begin{proposition}\label{H2estimate}
Let $u\in\CD_\Xi$ and $s>0$. For any $\delta>0$, we have
$$
(1-\delta)\|u_s^\sharp\|_{\CH^2}\le\|Hu\|_{L^2}+m_\delta^2(\Xi,s)\|u\|_{L^2}
$$
and
$$
\|Hu\|_{L^2}\le(1+\delta)\|u_s^\sharp\|_{\CH^2}+m_\delta^2(\Xi,s)\|u\|_{L^2}
$$
with $u_s^\sharp=\Phi^s(u)$ and $m_\delta^2(\Xi,s)>0$ an explicit constant.
\end{proposition}

\medskip

In addition to this comparison between $H$ and $L$ in norm, one has a similar statement in the quadratic form setting.

\medskip

\begin{proposition}\label{prop:formdom}
Let $u\in\CD_\Xi$ and $s>0$. For any $\delta>0$, we have
$$
(1-\delta)\langle\nabla u_s^\sharp,\nabla u_s^\sharp\rangle\le\langle u,Hu\rangle+m_\delta^1(\Xi,s)\|u\|_{L^2}^2
$$
and
$$
\langle u,Hu\rangle\le(1+\delta)\langle\nabla u_s^\sharp,\nabla u_s^\sharp\rangle+m_\delta^1(\Xi,s)\|u\|_{L^2}^2
$$
where $u_s^\sharp=\Phi^s(u)$ and $m_\delta^1(\Xi,s)>0$ an explicit constant. 
\end{proposition}

\medskip

One can show that $H\Gamma$ is the limit in norm of $H_\varepsilon\Gamma_\varepsilon$ as operators from $\CH^2$ to $L^2$ where
$$
H_\varepsilon:=L+\xi_\varepsilon-c_\varepsilon
$$
with $c_\varepsilon$ a diverging function as $\varepsilon$ goes to $0$, again see Section 2.1 of \cite{Mouzard}. In particular, one can take shift $c_\varepsilon$ by a large enough constant to ensure that $H$ is positive. Thus the previous Proposition implies that $\|\sqrt{H}u\|_{L^2}$ and $\|u_s^\sharp\|_{\CH^1}$ are equivalent. The diverging quantity is needed to take care of the singularity as explained in the introduction, this is the renormalisation procedure with
$$
\PI(X_1,\xi):=\lim_{\varepsilon\to0}\PI(X_1^{(\varepsilon)},\xi_\varepsilon)-c_\varepsilon
$$
in $\CC^{2\alpha-2}$. In the case of the torus, the noise is invariant by translation and the function $c_\varepsilon$ is actually a constant that diverges as $|\log\varepsilon|$, see \cite{AllezChouk}. This allows to prove that $H$ is a symmetric operator as the weak limit of the symmetric operators $H_\varepsilon$. Being closed and symmetric, it is enough to prove that
$$
(H+k)u=v
$$
admits a solution for some $k\in\IR$ to get self-adjointness for $H$, see Theorem $X.1$ in \cite{RS2}. This is done using the Babuška-Lax-Milgram Theorem, see \cite{Babuska} and Proposition \ref{prop:formdom} which implies that $H$ is almost surely bounded below. This implies self-adjointness and since the resolvent is a compact operator from $L^2$ to itself since $\CD_\Xi\subset\CH^\beta$ for any $\beta\in[0,\alpha)$.

\medskip

\begin{corollary}\label{SpectralResult}
The operator $H$ is self-adjoint with discrete spectrum $\big(\lambda_n(\Xi)\big)_{n\ge1}$ which is a nondecreasing diverging sequence without accumulation points. Moreover, we have
$$
L^2=\underset{n\ge1}{\bigoplus}\ \textup{Ker}\big(H-\lambda_n(\Xi)\big)
$$
with each kernel being of finite dimension. We finally have the min-max principle
$$
\lambda_n(\Xi)=\inf_D\sup_{u\in D;\|u\|_{L^2}=1}\langle Hu,u\rangle
$$
where $D$ is any $n$-dimensional subspace of $\CD_\Xi$; this can also be written as
$$
\lambda_n(\Xi)=\sup_{v_1,\ldots,v_{n-1}\in L^2}\ \inf_{\underset{\|u\|_{L^2}=1}{u\in\textup{Vect}(v_1,\ldots,v_{n-1})^\bot}}\langle Hu,u\rangle.
$$
\end{corollary}

\medskip

While the regularity of a function can be measured by its coefficients in the basis of eigenfunction of the Laplacian, the same is true for the Anderson Hamiltonian and the spaces agree if the regularity one considers is below the form domain.

\medskip

\begin{proposition}\label{HSobolevBound}
For $\beta\in(-\alpha,\alpha)$, there exists two constants $c_\Xi,C_\Xi>0$ such that
$$
c_\Xi\|H^{\frac{\beta}{2}}u\|_{L^2}\le\|u\|_{\CH^\beta}\le C_\Xi\|H^{\frac{\beta}{2}}u\|_{L^2}.
$$	
\end{proposition}

\medskip

\begin{proof}
Observe first that the statement is clear for $\beta=0,$ we consider only the case $\beta\in(0,\alpha)$ since the case of negative $\beta$ follows by duality. Again we take $(\varphi_n)_{n\ge1}$ and $(e_n)_{n\ge1}$ to denote the basis of eigenfunctions of $-\Delta$ and $H$ respectively. We have for any $v\in\CD_\Xi$
\begin{align*}
\big\|H^{\frac{\beta}{2}}v\big\|_{L^2}&=\Big(\sum_{n\ge1}\lambda_n^\beta\langle v,e_n\rangle^2\Big)^{\frac{1}{2}}\\
&=\Big(\sum_{n\ge1}\lambda_n^\beta\langle v,e_n\rangle^{2\beta}\langle v,e_n\rangle^{2-2\beta}\Big)^{\frac{1}{2}}\\
&\lesssim\big(\sum_{n\ge1}\lambda_n\langle v,e_n\rangle^2\big)^{\frac{\beta}{2}}\big(\sum_{n\ge1}\langle v,e_n\rangle^2\big)^{\frac{1-\beta}{2}}\\
&\lesssim\|H^{\frac{1}{2}}v\|_{L^2}^\beta\|v\|_{L^2}^{1-\beta}
\end{align*}
using Hölder's inequality. Thus the equivalence of $\|H^{\frac{1}{2}}v\|_{L^2}$  and $\|v_s^\sharp\|_{\CH^1}$ from Proposition \ref{prop:formdom}, together with the continuity of $\Phi^s$ from $L^2$ to itself yields
$$
\big\|H^{\frac{\beta}{2}}v\big\|_{L^2}\lesssim\|v_s^\sharp\|_{\CH^1}^\beta\|v^\sharp\|_{L^2}^{1-\beta}.
$$
Applying this with $v=\Gamma\big(\langle u_s^\sharp,\varphi_n\rangle\varphi_n\big)$ gives
\begin{align*}
\big\|H^{\frac{\beta}{2}}\Gamma\big(\langle u^\sharp,\varphi_n\rangle\varphi_n\big)\big\|_{L^2}&\lesssim\|\langle u^\sharp,\varphi_n\rangle\varphi_n\|_{\CH^1}^\beta\|\langle u^\sharp,\varphi_n\rangle\varphi_n\|_{L^2}^{1-\beta}\\
&\lesssim|\langle u^\sharp,\varphi_n\rangle|\|\varphi_n\|_{\CH^\beta}
\end{align*}
Thus
\begin{align*}
\|H^{\frac{\beta}{2}}u\|_{L^2}^2=\|H^{\frac{\beta}{2}}\Gamma(u_s^\sharp)\|_{L^2}^2&\le\sum_{n\ge1}\|H^{\frac{\beta}{2}}\Gamma\big(\langle u_s^\sharp,\varphi_n\rangle\varphi_n\big)\|_{L^2}^2\\
&\lesssim\sum_{n\ge1}|\langle u_s^\sharp,\varphi_n\rangle|^2\|\varphi_n\|_{\CH^{\beta}}^2\\
&\lesssim\|u_s^\sharp\|_{\CH^\beta}^2.
\end{align*}
Since $\beta\in[0,\alpha)$, we get
$$
\|H^{\frac{\beta}{2}}u\|_{L^2}\lesssim\|u\|_{\CH^\beta}.
$$
from the boundedness of $\Gamma$, see Proposition \ref{prop:gamma}.
The other inequality follows from the same reasoning with
$$
\|v\|_{\CH^\beta}\lesssim\|v_s^\sharp\|_{\CH^\beta}\lesssim\|v_s^\sharp\|_{\CH^1}^\beta\|v_s^\sharp\|_{L^2}^{1-\beta}\lesssim\|H^{\frac{1}{2}}v\|_{\CH^1}^\beta\|u\|_{L^2}^{1-\beta}
$$
and applying this bound to $u=\sum_{n\ge 1} \langle u,e_n \rangle e_n$ and proceeding as above we get the other direction.
\end{proof}

\medskip

The operator $H$ and its spectrum do not depend on $s>0$ but the different representation of $H$ as
$$
Hu=L\Phi^s(u)+\P_\xi\Phi^s(u)+\PI\big(\Phi^s(u),\xi\big)+R(u)+\Psi^s(u)
$$
yields different bounds on the eigenvalues. We state the simpler form for the bounds, see \cite{Mouzard} for the general result. It is sharp enough to obtain an almost sure Weyl-type law from the one for the Laplace-Beltrami operator.

\medskip

\begin{proposition}\label{EingenBoundWeylLaw}
Let $\delta\in(0,1)$. Then there exists two constants $m_\delta^1(\Xi),m_\delta^2(\Xi)$ such that
$$
\lambda_n-m_\delta^1(\Xi)\le\lambda_n(\Xi)\le(1+\delta)\lambda_n+m_\delta^2(\Xi)
$$
for any $n\in\IN$. This implies the almost sure Weyl-type law
$$
\lim_{\lambda\to\infty}\lambda^{-1}|\{n\ge0;\lambda_n(\Xi)\le\lambda\}|=\frac{\textup{Vol}(M)}{4\pi}.
$$
\end{proposition}

\bigskip

\section{Strichartz inequalities for the stochastic Schrödinger equation}\label{SectionStrichartzSchr}

For the rest of the work, we fix a parameter $s>0$ small enough in order to have all the needed continuity estimates. Every constant may implicitly depend on $s$ and on the norm of the enhanced noise, we do not explicate the dependence since it is not relevant at this stage. From now on, we will also use that $\alpha$ can be taken arbitrary close to $1$ since it is given by the regularity of the spatial white noise. We consider the Schrödinger operator
$$
\G:=\Gamma^{-1}H\Gamma
$$
which appears naturally when transforming the Schrödinger equation and the wave equation with multiplicative noise. In fact, if $u$ solves
$$
\left|\begin{array}{ccc}
i\partial_tu+Hu&=&0,\\
u(0)&=&u_0
\end{array}\right.
$$
then $u^\sharp:= \Gamma^{-1}u$ solves the transformed equation
$$
\left|\begin{array}{ccc}
i\partial_tu^\sharp+H^\sharp u^\sharp&=&0,\\
u^\sharp(0)&=&\Gamma^{-1} u_0
\end{array}\right.
$$
In this Section, we show Strichartz inequalities for the associated Schrödinger equation with an arbitrary small loss of regularity with respect to the deterministic case. Afterwards, in Section \ref{subsec:wellNLS}, we detail how these can be used to get a low-regularity solution theory for the nonlinear Schrödinger equation with multiplicative noise.

\bigskip

\subsection{Strichartz inequalities for the Schrödinger group}

As was hinted at in Proposition \ref{H2estimate}, the transformed operator $H^\sharp$ it is a lower-order perturbation of the Laplace-Beltrami operator. We obtain the following result which is somewhat similar to Theorem $6$ in the work \cite{BGT} by Burq, Gérard and Tzvetkov, where they proved that the Strichartz inequalities are stable for some lower order perturbations. This does not cover the case of the Anderson Hamiltonian however our proof is very similar, see also \cite{Zachhuber}.

\medskip

\begin{proposition}\label{prop:opdiff}
Let $0\le\beta<1$. For any $\kappa>0$, we have
$$
\|(\G-L)v\|_{\CH^\beta}\lesssim\|v\|_{\CH^{1+\beta+\kappa}}.
$$
\end{proposition}

\medskip

\begin{proof}
For $u=\Gamma u^\sharp\in\CD_\Xi$, recall that
$$
Hu=Lu^\sharp+\P_\xi u^\sharp+\PI(u^\sharp,\xi)+R(u)
$$
where
\begin{align*}
R(u)&:=\PI\big(u,\PI(X_1,\xi)\big)+\P_{\PI(X_1,\xi)}u+\DC(u,X_1,\xi)+\P_u\PI(X_2,\xi)+\DD(u,X_2,\xi)\\
&\quad+\DS(u,X_2,\xi)+\P_\xi\PT_uX_2-e^{-L}\left(\P_uX_1+\P_uX_2\right).
\end{align*}
Thus $\G v$ is given by
$$
\G v=Lv+\P_\xi v+\PI(v,\xi)+R(\Gamma v)-\PT_{H\Gamma v}(X_1+X_2)
$$
and for any $\kappa>0$ and $\beta\in[0,\alpha]$, we have
\begin{align*}
\|(\G-L)v\|_{\CH^\beta}&\lesssim\|\P_\xi v+\PI(v,\xi)\|_{\CH^\beta}+\|R(\Gamma v)\|_{\CH^\beta}+\|\PT_{H\Gamma v}(X_1+X_2)\|_{\CH^\beta}\\
&\lesssim\|\xi\|_{\CC^{-1-\kappa}}\|v\|_{\CC^{\beta+1+\kappa}}+\|\Gamma v\|_{\CH^\alpha}+\|H\Gamma v\|_{\CH^{-1+\kappa+\beta}}\|X_1+X_2\|_{\CH^{1-\kappa}}\\
&\lesssim\|v\|_{\CH^{1+\beta+\kappa}}+\|v\|_{\CH^\alpha}+\|v\|_{\CH^{1+\kappa+\beta}}
\end{align*}
using Proposition \ref{SobolevBoundH} and the proof is complete since $\alpha<1$.
\end{proof}

\medskip

Since the unitary group associated to $H$ is bounded on $L^2$ and on the domain $\CD_\Xi$ of $H$, this implies a similar result for the ``sharpened'' group associated with $\G$ in terms of classical Sobolev spaces. Recall that $\G=\Gamma^{-1}H\Gamma$ with $\Gamma$ an isomorphism from $L^2$ to itself thus $e^{it\G}:=\Gamma^{-1}e^{itH}\Gamma$ is well-defined on $L^2$. We now state some of its properties.

\medskip

\begin{proposition}\label{prop:sharpgroup}
For any $0\le\beta\le 2$ and $t\in\IR$, we have
$$
\|e^{it\G}v\|_{\CH^\beta}\lesssim\|v\|_{\CH^\beta}.
$$
Moreover, $e^{itH^\sharp}$ is a non-unitary strongly continuous group of $L^2$ bounded operators, namely
$$
e^{i(t+s)H^\sharp}v=e^{itH^\sharp}e^{isH^\sharp}v
$$
for all $s,t\in\mathbb{R}$ and $v\in L^2$.
\end{proposition}

\medskip

\begin{proof}
For $\beta=0$, this follows from the continuity of $\Gamma$ and $\Gamma^{-1}$ from $L^2$ to itself. For $\beta=2$, we have
\begin{align*}
\|e^{it\G}v\|_{\CH^2}&=\|\Gamma^{-1}e^{itH}\Gamma v\|_{\CH^2}\\
&\lesssim\|He^{itH}\Gamma v\|_{L^2}\\
&\lesssim\|e^{itH}H\Gamma v\|_{L^2}\\
&\lesssim\|H\Gamma v\|_{L^2}\\
&\lesssim\|v\|_{\CH^2},
\end{align*}
having used Proposition \ref{H2estimate}. The results for any $\beta\in(0,2)$ is obtained by interpolation and the group property follows simply from the group property of $e^{itH}$ by observing
$$
e^{i(t+s)H^\sharp} v=\Gamma^{-1}e^{i(t+s)H}\Gamma v=\Gamma^{-1}e^{itH}\Gamma\Gamma^{-1}e^{isH}\Gamma v=e^{itH^\sharp} e^{isH^\sharp v}
$$
\end{proof}

\medskip

Strichartz inequalities are refinements of the estimates from the previous Proposition. The following statement is such a result, which has an arbitraty small loss of derivative coming from the irregularity of the noise in the addition to the $\frac{1}{p}$ loss from the manifold setting without boundary which one sees in \cite{BGT}. We refer to a pair $(p,q)$ statisfying
$$
\frac{2}{p}+\frac{2}{q}=1
$$
as a Strichartz pair from now on.

\medskip

\begin{theorem}\label{StrichartzSch}
Let $M$ be a two-dimensional compact manifold without boundary and let $(p,q)$ be a Strichartz pair. Then for any $\varepsilon>0$
$$
\|e^{it\G}v\|_{L^p([0,1],L^q)}\lesssim\|v\|_{\CH^{\frac{1}{p}+\varepsilon}}.
$$
This implies the bound
$$
\|e^{itH}u\|_{L^p([0,1],L^q)}\lesssim\|\Gamma^{-1} u\|_{\CH^{\frac{1}{p}+\varepsilon}}\lesssim\|u\|_{\CH^{\frac{1}{p}+\varepsilon}}.
$$
\end{theorem}

\medskip

First, we need to prove the following Lemma. It gives the difference between the Schrödinger groups associated to $\G$ and $L$ from the difference between $\G$ and $L$ itself. Moreover it quantifies that their difference is small in a short time interval if one gives up some regularity.

\medskip

\begin{lemma}\label{GroupDiff}
Given $v\in\CH^2$, we have
$$
\Big(e^{i(t-t_0)\G}-e^{i(t-t_0)L}\Big)v=i\int_{t_0}^te^{i(t-s)L}(\G-L)e^{i(s-t_0)\G}v\drm s
$$
for any $t,t_0\in\IR$.
\end{lemma}

\medskip

\begin{proof}
The ``sharpened'' group yields the solution of the Schrödinger equation
$$
\big(i\partial_t+\G\big)(e^{i(t-t_0)\G}v)=0
$$
which is equal to $v$ at time $t=t_0$ thus
$$
\big(i\partial_t+L\big)(e^{i(t-t_0)\G}v)=\big(L-\G\big)(e^{i(t-t_0)\G}v).
$$
Using the unitary group representation of the solution to the Schrödinger equation associated to $L$, we deduce that
$$
\big(i\partial_t+L\big)(e^{i(t-t_0)L}v-e^{i(t-t_0)\G}v)=\big(\G-L\big)(e^{i(t-t_0)\G}v).
$$
Since the solution is equal to $0$ at time $t=t_0$, the mild formulation of this last equation yields
$$
\Big(e^{i(t-t_0)\G}-e^{i(t-t_0)L}\Big)v=i\int_{t_0}^te^{i(t-s)L}(\G-L)e^{i(s-t_0)\G}v\drm s.
$$
\end{proof}

\medskip

\begin{proofof}{Theorem}{\ref{StrichartzSch}}
For $N\in\IN^*$ to be chosen, we have
$$
\|e^{it\G}v\|_{L^p([0,1],L^q)}^p=\sum_{\ell=0}^N\|e^{it\G}v\|_{L^p([t_\ell,t_{\ell+1}],L^q)}^p
$$
where $t_\ell:=\frac{\ell}{N}$. For $t\in[t_\ell,t_{\ell+1})$, the previous Lemma gives
$$
e^{it\G}v=e^{i(t-t_\ell)\G}e^{it_\ell\G}v=e^{i(t-t_\ell)L}e^{it_\ell\G}v+i\int_{t_\ell}^te^{i(t-s)L}(\G-L)e^{is\G}v\drm s.
$$
Applying this with $v=\Delta_ku$ gives
\begin{align*}
\|\Delta_je^{it\G}\Delta_ku\|_{L^p([0,1],L^q)}^p&\le\sum_{\ell=0}^N\big\|\Delta_je^{i(t-t_\ell)L}e^{it_\ell\G}\Delta_ku\big\|_{L^p([t_\ell,t_{\ell+1}],L^q)}^p\\
&\quad+\sum_{\ell=0}^N\Big\|\Delta_j\int_{t_\ell}^te^{i(t-s)L}(\G-L)e^{is\G}\Delta_ku\drm s\Big\|_{L^p([t_\ell,t_{\ell+1}],L^q)}^p.
\end{align*}
Assume $N\ge 2^j$ such that $|t_{\ell+1}-t_\ell|\le2^{-j}$. For the first term, we have
\begin{align*}
\|\Delta_je^{i(t-t_\ell)L}e^{it_\ell\G}\Delta_ku\|_{L^p([t_\ell,t_{\ell+1}],L^q)}^p&=\|e^{i(t-t_\ell)L}\Delta_je^{it_\ell\G}\Delta_ku\|_{L^p([t_\ell,t_{\ell+1}],L^q)}^p\\
&\lesssim\|\Delta_je^{it_\ell\G}\Delta_ku\|_{L^2}^p\\
&\lesssim2^{-j\delta p}\|\Delta_je^{it_\ell\G}\Delta_ku\|_{\CH^\delta}^p\\
&\lesssim2^{-j\delta p}2^{-k\delta'p}\|\Delta_ku\|_{\CH^{\delta+\delta'}}^p
\end{align*}
for any $\delta,\delta'\in\IR$ using Proposition \ref{prop:sharpgroup}, Strichartz inequality for spectrally localised data from Section \ref{SubsectionStrichartz} and Bernstein's Lemma, see Lemma \ref{lem:bern}. For the second term, we have
\begin{align*}
\Big\|\Delta_j\int_{t_\ell}^te^{i(t-s)L}(\G-L)e^{is\G}\Delta_ku\drm s\Big\|_{L^p([t_\ell,t_{\ell+1}],L^q)}^p&=\Big\|\int_{t_\ell}^te^{i(t-s)L}\Delta_j(\G-L)e^{is\G}\Delta_ku\drm s\Big\|_{L^p([t_\ell,t_{\ell+1}],L^q)}^p\\
&\lesssim\left(\int_{t_\ell}^{t_{\ell+1}}\big\|e^{i(t-s)L}\Delta_j(\G-L)e^{is\G}\Delta_ku\big\|_{L^p([t_\ell,t_{\ell+1}],L^q)}\drm s\right)^p\\
&\lesssim\left(\int_{t_\ell}^{t_{\ell+1}}\big\|\Delta_j(\G-L)e^{is\G}\Delta_ku\big\|_{L^2}\drm s\right)^p\\
&\lesssim2^{-j\sigma p}\left(\int_{t_\ell}^{t_{\ell+1}}\big\|\Delta_j(\G-L)e^{is\G}\Delta_ku\big\|_{\CH^\sigma}\drm s\right)^p\\
&\lesssim2^{-j\sigma p}\left(\int_{t_\ell}^{t_{\ell+1}}\big\|(\G-L)e^{is\G}\Delta_ku\big\|_{\CH^\sigma}\drm s\right)^p\\
&\lesssim2^{-j\sigma p}\left(\int_{t_\ell}^{t_{\ell+1}}\big\|e^{is\G}\Delta_ku\big\|_{\CH^{\sigma+1+\kappa}}\drm s\right)^p\\
&\lesssim N^{-p}2^{-j\sigma p}\|\Delta_ku\|_{\CH^{1+\sigma+\kappa}}^p\\
&\lesssim N^{-p}2^{-j\sigma p}2^{-k\sigma'p}\|\Delta_ku\|_{\CH^{1+\sigma+\sigma'+\kappa}}^p
\end{align*}
for any $\sigma\in(0,1),\sigma'\in\IR$ and $0<\kappa<1-\alpha$ where again the dyadic factors come from Bernstein's Lemma and we have used the bounds from Propositions \ref{prop:opdiff} and \ref{prop:sharpgroup} with Strichartz inequality for spectrally localised data. Summing over the sub-intervals gives
$$
\|\Delta_je^{it\G}\Delta_ku\|_{L^p([0,1],L^q)}\lesssim N^{\frac{1}{p}}2^{-j\delta}2^{-k\delta'}\|\Delta_ku\|_{\CH^{\delta+\delta'}}+N^{\frac{1-p}{p}}2^{-j\sigma}2^{-k\sigma'}\|\Delta_ku\|_{\CH^{1+\sigma+\sigma'+\kappa}}.
$$
Let $\eta>0$ small and take 
$$
N=2^j,\ \delta=\eta+\frac{1}{p},\ \delta'=\sigma=\sigma'=\eta
$$
which satisfies in particular $N\ge 2^j$ and $\sigma\in[0,\alpha)$ to sum over $k\le j$. We get
\begin{align*}
\|\sum_{k\le j}\Delta_je^{it\G}\Delta_ku\|_{L^p([0,1],L^q)}&\lesssim \sum_j \sum_{k\le j}\|\Delta_je^{it\G}\Delta_ku\|_{L^p([0,1],L^q)}\\
&\lesssim \sum_{j\ge 0} \sum_{k\le j} 2^{\frac{j}{p}}2^{-j(\frac{1}{p}+\eta)}2^{-k\eta}\|\Delta_ku\|_{\CH^{\frac{1}{p}+2\eta}}+2^{j\frac{1-p}{p}}2^{-j\eta}2^{-k\eta}\|\Delta_ku\|_{\CH^{1+2\eta+\kappa}}\\
&\lesssim \sum_{j\ge 0}  2^{-j\eta}\|\Delta_{\le j}u\|_{\CH^{\frac{1}{p}+2\eta}}+2^{j\frac{1-p}{p}}2^{-j\eta}\|\Delta_{\le j}u\|_{\CH^{1+2\eta+\kappa}}\\
&\lesssim \|u\|_{\CH^{\frac{1}{p}+2\eta}}+\sum_{j\ge 0}2^{-j\eta} 2^{j\frac{1-p}{p}}2^{-j\frac{1-p}{p}}\|\Delta_{\le j}u\|_{\CH^{1-1+\frac{1}{p}+2\eta+\kappa}}\\
&\lesssim \|u\|_{\CH^{\frac{1}{p}+2\eta}}+ \|u\|_{\CH^{\frac{1}{p}+2\eta+\kappa}},
\end{align*}
having used Bernstein's inequality, Lemma \ref{lem:bern}, for the projector $\Delta_{\le j}$. For the sum over $j\le k$, we choose instead 
$$
N=2^k,\ \delta= \delta'=\sigma=\sigma'=\eta,
$$
with $\eta>0$ small as before. Since $j\le k$, we have $N\ge2^j$ thus get the bound for the other part of the double sum
\begin{align*}
\|\sum_{j\le k}\Delta_je^{it\G}\Delta_ku\|_{L^p([0,1],L^q)}&\lesssim \sum_{k\ge 0} \sum_{j\le k}\|\Delta_je^{it\G}\Delta_ku\|_{L^p([0,1],L^q)}\\
&\lesssim \sum_{k\ge 0} \sum_{j\le k} 2^{\frac{k}{p}}2^{-j\eta}2^{-k\eta}\|\Delta_ku\|_{\CH^{2\eta}}+2^{\frac{k(1-p)}{p}}2^{-j\eta}2^{-k\eta}\|\Delta_ku\|_{\CH^{1+2\eta+\kappa}}\\
&\lesssim \|u\|_{\CH^{\frac{1}{p}+2\eta}}+\|u\|_{\CH^{1+\frac{1-p}{p}+2\eta+\kappa}}\\
&\lesssim \|u\|_{\CH^{\frac{1}{p}+2\eta+\kappa}},
\end{align*}
having used Bernstein's inequality again. This completes the proof since $\eta$ and $\kappa$ can be taken arbitrary small. This implies the bound
$$
\|e^{itH}u\|_{L^p([0,1],L^q)}\lesssim\|\Gamma^{-1}u\|_{\CH^{\frac{1}{p}+\varepsilon}}\lesssim\|u\|_{\CH^{\frac{1}{p}+\varepsilon}}
$$
using that $\frac{1}{p}+\varepsilon<1$ and Proposition \ref{prop:gamma}.
\end{proofof}

\medskip

\begin{remark}
We proved that Strichartz inequalities are stable under suitable perturbation, that is lower-order perturbation in the sense of previous Proposition \ref{prop:opdiff}. This is similar in spirit to Theorem $6$ in \cite{BGT}. One can show that the magnetic Laplacian with white noise magnetic field constructed in \cite{MM} is also a lower order perturbation of the Laplacian on the two-dimensional torus in this sense. Thus Theorem \ref{StrichartzSch} also gives Strichartz inequalities for the associated Schrödinger group.
\end{remark}

\medskip

As a corollary, we state the inhomogeneous inequalities needed to solve the nonlinear equation. This is straightforward, see \cite{Zachhuber} and references therein.

\medskip

\begin{corollary}\label{cor:inhom}
In the setting of Theorem \ref{StrichartzSch}, we have in addition the bound
\begin{align*}
\left\|\int_0^te^{i(t-s)\G}f(s)\drm s\right\|_{L^p([0,1],L^q)}\lesssim\int_0^1\|f(s)\|_{\CH^{\frac{1}{p}+\varepsilon}}\drm s
\end{align*} 
for all $f\in L^1([0,1],\mathcal{H}^{\frac{1}{p}+\varepsilon})$.
\end{corollary}

\medskip

The only ingredient in the proof of the Theorem where the boundary appears is when we apply the result for the Laplacian. By using Theorem \ref{thm:Strichartzschboundary} and in particular \eqref{Strichartzloc2} instead of \eqref{Strichartzloc}, we immediately get the following analogous result which is of course weaker.

\medskip

\begin{theorem}\label{thm: strboundary}
Let $M$ be a compact surface with boundary. Let $p\in(3,\infty]$ and $q\in[2,\infty)$ such that
$$
\frac{3}{p}+\frac{2}{q}=1.
$$
Then for any $\varepsilon>0$
$$
\|e^{it\G}u\|_{L^p([0,1],L^q)}\lesssim\|u\|_{\CH^{\frac{2}{p}+\varepsilon}}.
$$
and 
\begin{align*}
\left\|\int_0^te^{i(t-s)\G}f(s)\drm s\right\|_{L^p([0,1],L^q)}\lesssim\int_0^1\|f(s)\|_{\CH^{\frac{2}{p}+\varepsilon}}\drm s
\end{align*} 
for all $f\in L^1([0,1],\mathcal{H}^{\frac{2}{p}+\varepsilon})$.
\end{theorem}

\bigskip

\subsection{Local well-posedness for multiplicative stochastic cubic NLS}\label{subsec:wellNLS}

We now apply our results to the local in time well-posedness of the cubic multiplicative stochastic NLS
\begin{equation*}
\label{eqn:cubicNLS}
\left|\begin{array}{ccc}
i\partial_t u+H u&=&-|u|^2u\\
u(0)&=&u_0
\end{array}\right.
\end{equation*}
with $u_0\in\CH^\sigma$ where $\sigma\in(\frac{1}{2},1)$ and in the energy space, that is $u_0\in \mathcal{D}(\sqrt{H})=\Gamma \mathcal{H}^1$. The latter hypothesis is natural to assume, since solutions starting in the energy space, usually called energy solutions, are intimately related to the conserved energy
\begin{equation}
E(u)(t):=\frac{1}{2}(u(t),H u(t))+\frac{1}{4}\int|u(t)|^4=E(u_0),
\label{eqn:nlsenergy}
\end{equation}
introduced in \cite{GUZ} on the torus and \cite{Mouzard} on general surfaces.  Thus we refer to $\CD(\sqrt{H})=\Gamma\CH^1$ as the energy space for the Anderson Hamiltonian, see Proposition \ref{prop:formdom}. Note that, as is usual in these types of fixed point arguments, the sign of the nonlinearity does not play a role for local well-posedness, but for the sake of definiteness we take the defocussing nonlinearity. We remark also that one can prove similar results for more general nonlinearities, we considered only the cubic equation in this work. See for example \cite{TV,TVbis} where they considered generic polynomial nonlinearity and obtain global well-posedness. As explained in Section $3.2.2$ of \cite{GUZ}, their result for the equation with white noise potential is weaker than the one for the deterministic equation since Strichartz inequalities were not know in this singular case. This was a motivation for the study of Strichartz estimates for the operator $H$. The well-posedness itself follows from a fairly straightforward contraction argument, similar to e.g. Proposition 3.1 in \cite{BGT}. Finally, we only consider a surface without boundary, the case with boundary is analogous using Theorem \ref{thm: strboundary} instead of Theorem \ref{StrichartzSch}. The mild formulation is
$$\label{eqn:mildH}
u(t)=e^{itH}u_0-i\int_{0}^{t}e^{i(t-s)H}\big(|u|^2u\big)(s)\drm s
$$
and applying the $\Gamma^{-1}$ map introduced in Section \ref{SubsectionAnderson} yields the mild formulation for the transformed quantity $u^\sharp=\Gamma^{-1}u$. We get
$$\label{eqn:mildsharp}
u^\sharp(t)=e^{itH^\sharp}u^\sharp_0-i\int_{0}^{t}e^{i(t-s)H^\sharp}\Gamma^{-1}\big(|\Gamma u^\sharp|^2\Gamma u^\sharp\big)(s)\drm s
$$
where $u_0^\sharp:=\Gamma^{-1}u_0$, this is where the transformed operator $H^\sharp=\Gamma^{-1}H\Gamma$ appears naturally. Despite the seemingly complicated nonlinear expression, this new mild formulation is easier to deal with since $H^\sharp$ is a perturbation of the Laplacian and has domain $\mathcal{H}^2,$ hence it is not as outlandish as $H$ and its domain which contains no non-zero smooth functions. Now, we have to find a bound for the map
$$
\Psi (v)(t):=e^{itH^\sharp}v_0-i\int_{0}^{t}e^{i(t-s)H^\sharp}\Gamma^{-1}\big(|\Gamma v|^2\Gamma v\big)(s)\drm s
$$
in a suitable space which allows us to get a unique fixed point. One then recovers a solution to the original equation with $u:=\Gamma v$ and choosing $v_0:=\Gamma^{-1}u_0$. Since $\Gamma$ is an isomorphism on $L^p$ for $p\in[2,\infty]$ and both $\Gamma$ and $e^{it\G}$ are isomorphisms on $\CH^\sigma$ to itself for $\sigma\in[0,1)$, it is natural to consider initial datum $v_0$, and thus also $u_0$, in $\mathcal{H}^\sigma$ for $0<\sigma<1.$ Therefore we bound $\Psi(v)$ in $\mathcal{H}^\sigma$ with
\begin{align*}
\big\|\Psi(v)(t)\big\|_{\mathcal{H}^\sigma}&\lesssim \|v_0\|_{\mathcal{H}^\sigma}+\int_{0}^{t}\|\Gamma v(s)^3\|_{\mathcal{H}^\sigma}\drm s\\
&\lesssim \|v_0\|_{\mathcal{H}^\sigma}+\int_{0}^{t}\|\Gamma v(s)\|_{\mathcal{H}^\sigma} \|\Gamma v(s)\|^2_{L^\infty}\drm s\\
&\lesssim \|v_0\|_{\mathcal{H}^\sigma}+\int_{0}^{t}\|v(s)\|_{\mathcal{H}^\sigma}\| v(s)\|^2_{L^\infty}\drm s\\
&\lesssim \|v_0\|_{\mathcal{H}^\sigma}+\|v\|_{L^\infty([0,t],\mathcal{H}^\sigma)}\|v\|^2_{L^2([0,t],L^\infty)}
\end{align*}
where in the first and third lines we have used the continuity of $e^{it\G}$ and $\Gamma$ and Lemma \ref{lem:tame} in the second line. For $\sigma<1$, the space $\CH^\sigma$ is not an algebra and one can not simply use its norm to bound the nonlinearity. However, one may bound it using the $L^\infty$-norm in space by observing that one needs less integrability in time and this is precisely the point where the Strichartz estimates turn out to be useful. As for the deterministic equation, we work with the function spaces
$$
\CW^{\beta,q}(M)=\{u\in\CD'(M);(1-\Delta)^{\frac{\beta}{2}}u\in L^q\}
$$
with associated norm
$$
\|u\|_{\CW^{\beta,q}}:=\|(1-\Delta)^{\frac{\beta}{2}}u\|_{L^q}.
$$
For $\beta\in[0,1)$ and $q=2$, one recovers the Sobolev spaces and the norm is equivalent to
$$
\|u\|_{\CH^\beta}=\|(1+H)^{\frac{\beta}{2}}u\|_{L^2}
$$
by Proposition \ref{HSobolevBound}. Within this framework, Strichartz inequalities from Theorem \ref{StrichartzSch} gives us the bound
$$
\|e^{itH^\sharp} w\|_{L^p([0,1],\CW^{\beta,q})}\lesssim \|w\|_{\mathcal{H}^{\frac{1}{p}+\beta+\kappa}},
$$
for any Strichartz pair $(p,q)$ and $\kappa>0$. Furthermore, the space $\CW^{\beta,q}$ is continuously embedded in $L^\infty$ for $\beta>\frac{2}{q}$. Let $\sigma\in\IR$ such that
$$
\frac{1}{p}+\frac{2}{q}+2\kappa\le\sigma.
$$
Thus for $0<t\le 1$, we get the bound
\begin{align*}
\|\Psi(v)\|_{L^p\big([0,t],\CW^{\frac{2}{q}+\kappa,q}\big)}&\lesssim \|v_0\|_{\mathcal{H}^{\frac{1}{p}+\frac{2}{q}+2\kappa}}+\int_0^t\big\|\Gamma^{-1}\big(|\Gamma v|^2\Gamma v\big)(s)\big\|_{\mathcal{H}^{\frac{1}{p}+\frac{2}{q}+2\kappa}}\drm s\\
&\lesssim \|v_0\|_{\mathcal{H}^{\sigma}}+\int_0^t\big\|\Gamma v(s)^3\big\|_{\mathcal{H}^{\sigma}}\drm s\\
&\lesssim\|v_0\|_{\mathcal{H}^\sigma}+\|v\|_{L^\infty([0,t],\mathcal{H}^\sigma)}\|v\|^2_{L^2([0,t],L^\infty)}\\
&\lesssim \|v_0\|_{\mathcal{H}^\sigma}+\|v\|_{L^\infty([0,t],\mathcal{H}^\sigma)}\|v\|^2_{L^2\big([0,t],\CW^{\frac{2}{q}+\kappa,q}\big)}\\
&\lesssim\|v_0\|_{\mathcal{H}^\sigma}+t^{\frac{p-2}{p}}\|v\|_{L^\infty([0,t],\mathcal{H}^\sigma)}\|v\|^2_{L^p\big([0,t],\CW^{\frac{2}{q}+\kappa,q}\big)}
\end{align*}
using Corollary \ref{cor:inhom} in the first line, Hölder inequality in the last line and bicontinuity of $\Gamma$ from $\CH^\sigma$ to itself. For $0<t'\le t$, we also have
\begin{align*}
\big\|\Psi(v)(t')\big\|_{\mathcal{H}^\sigma}&\lesssim \|v_0\|_{\mathcal{H}^\sigma}+\int_0^{t'}\|v(s)\|_{\mathcal{H}^\sigma}\|v(s)\|^2_{L^\infty}\drm s\\
&\lesssim \|v_0\|_{\mathcal{H}^\sigma}+\|v\|_{L^\infty([0,t'],\mathcal{H}^\sigma)}\|v\|^2_{L^2([0,t'],L^\infty)},\\
&\lesssim\|v_0\|_{\mathcal{H}^\sigma}+{t'}^{\frac{p-2}{p}}\|v\|_{L^\infty([0,t'],\mathcal{H}^\sigma)}\|v\|^2_{L^p\big([0,t'],\CW^{\frac{2}{q}+\kappa,q}\big)}.
\end{align*}
This gives us the combined bound
$$\label{eqn:boundSstr}
\|\Psi(v)\|_{L^p\big([0,t],\CW^{\frac{2}{q}+\kappa,q}\big)}+\|\Psi(v)\|_{L^\infty([0,t],\mathcal{H}^\sigma)}\lesssim \|v_0\|_{\mathcal{H}^\sigma}+t^{\frac{p-2}{p}}\|v\|_{L^\infty([0,t],\mathcal{H}^\sigma)}\|v\|^2_{L^p\big([0,t],\CW^{\frac{2}{q}+\kappa,q}\big)}
$$
that will be the main tool for the fixed point. Note that the restrictions
$$
\frac{1}{p}+\frac{2}{q}+2\kappa\le\sigma\quad\text{and}\quad\frac{2}{p}+\frac{2}{q}=1
$$
gives
$$
1-\frac{1}{p}+2\kappa\le\sigma.
$$
Since $p\ge2$ and $\kappa>0$ can be taken arbitrary small, this gives
$$
\sigma>\frac{1}{2}
$$
and leads to the following local-in-time well-posedness result. Without Strichartz estimates, even in the classical case, one could not go beyond the threshold $\sigma\ge1$.

\medskip

\begin{theorem}\label{thm:NLSlwp}
Let $M$ be a compact surface without boundary, $\sigma>\frac{1}{2}$ and initial data $v_0\in \mathcal{H}^\sigma$. Let $\kappa>0$ and $(p,q)$ a Strichartz pair such that
$$
\frac{1}{p}+\frac{2}{q}+2\kappa\le\sigma,
$$
there exists a time $T>0$ until which there exists a unique solution
$$
v\in C\big([0,T],\mathcal{H}^\sigma\big)\cap L^p\big([0,T],\CW^{\frac{2}{q}+\kappa,q}\big)
$$
to the mild formulation of the transformed PDE
$$
\left|\begin{array}{ccc}
i\partial_t v+H^\sharp v&=&-\Gamma^{-1}\big(|\Gamma v|^2\Gamma v\big)\\
v(0)&=&v_0
\end{array}\right..
$$
Moreover, the solution depends continuously on the initial data $v_0\in\CH^\sigma$.
\end{theorem}

\medskip

\begin{proof}
This is a straightforward contraction argument where the main ingredient is the bound proved in the preceding arguments. By choosing the radius $R$ of the ball and the final time appropriately, we can prove that
$$
\Psi:B(0,R)_{C([0,T],\mathcal{H}^\sigma)\cap L^p(([0,T],\CW^{\frac{2}{q}+\kappa,q})}\to B(0,R)_{C([0,T],\mathcal{H}^\sigma)\cap L^p(([0,T],\CW^{\frac{2}{q}+\kappa,q})}
$$
is in fact a contraction. Using the previously established bound \eqref{eqn:boundSstr} and an analogous bound for the difference, one finds that this can be achieved if one chooses
\begin{equation*}
R=2\tilde{C}\|v_0\|_{\CH^\sigma}\quad\text{and}\quad T=\left(\frac{1}{3R^2\tilde{C}}\right)^{\frac{p}{p-2}}
\end{equation*}
for some constant $\tilde{C}$ which depends on the norm of the enhanced noise $\Xi$ and the parameters appearing.
\end{proof}

\medskip

Finally we give the analogous result for surfaces with boundary which is of course weaker, however we still get a better result than one gets simply from using the algebra property of Sobolev spaces.

\medskip

\begin{theorem}\label{thm:NLSlwpbrdy}
Let $M$ be a compact surface with boundary, $\sigma>\frac{2}{3}$ and $p,q,\kappa$ s.t.
\begin{align*}
\frac{3}{p}+\frac{2}{q}=1 \text{ and }\frac{2}{p}+\frac{2}{q}+2\kappa\le \sigma.
\end{align*}
For any initial datum $v_0\in \mathcal{H}^\sigma$ there exists a unique solution 
\begin{align*}
v\in C([0,T],\mathcal{H}^\sigma)\cap L^p([0,T],\mathcal{W}^{\frac{2}{q}+\kappa,q})
\end{align*}
to the mild formulation of the transformed PDE up to a time $T>0$ depending on the data which depends continuously on the initial condition.
\end{theorem}

\bigskip

\section{Strichartz inequalities for the stochastic wave equation}

Again, we consider the ``sharpened'' operator
$$
\G:=\Gamma^{-1}H\Gamma
$$
which appears naturally when transforming the wave equation with multiplicative noise. If $u$ solves
$$
\left|\begin{array}{ccc}
\partial_t^2u+Hu&=&0\\
(u,\partial_tu)|_{t=0}&=&(u_0,u_1)
\end{array}\right.
$$
then $u^\sharp:= \Gamma^{-1}u$ solves the transformed equation
$$
\left|\begin{array}{ccc}
\partial_t^2u^\sharp+H^\sharp u^\sharp&=&0\\
(u^\sharp,\partial_tu^\sharp)|_{t=0}&=&(\Gamma^{-1}u_0,\Gamma^{-1}u_1)
\end{array}\right.
$$
In this Section, we show Strichartz inequalities for the associated wave equation. We will further detail how these can be used to get a low-regularity solution theory for the nonlinear wave equation with multiplicative noise. This equation was also considered in \cite{Zachhuberwave} by Zachhuber on the full space in two and three dimensions where global well-posedness is obtained using finite speed of propagation.

\bigskip

\subsection{Strichartz inequalities for the wave propagator}

The propagator associated to the wave equation is
$$
(u_0,u_1)\mapsto\cos(t\sqrt{H})u_0+\frac{\sin(t\sqrt{H})}{\sqrt{H}}u_1
$$
with initial conditons $(u,\partial_tu)|_{t=0}=(u_0,u_1)$. As for the Schrödinger equation, the following Strichartz inequalities hold on a two-dimensional compact Riemannian manifold without boundary, see \cite{BSSwave} and the references therein. We state the result in the homogeneous case for simplicity however one directly obtains inhomogenous bounds as in Corollary \ref{cor:inhom}. We cite the following Strichartz estimates which hold on compact surfaces respectively without and with boundary, see \cite{BSSwave}.

\medskip

\begin{theorem}
Let $(M,g)$ be a compact two-dimensional Riemannian manifold without boundary. Let $p,q\in[2,\infty]$ such that
$$
\frac{2}{p}+\frac{1}{q}\le\frac{1}{2}
$$
and consider
$$
\frac{1}{p}+\frac{2}{q}:=1-\sigma.
$$
Then the solution to
\begin{eqnarray*}
(\partial^2_t-\Delta_g)u & = & 0\\
(u,\partial_tu)|_{t=0} & = & (u_0,u_1)\in\mathcal{H}^{\sigma}\times\mathcal{H}^{\sigma-1}
\end{eqnarray*}
satisfies the bound
$$
\|u\|_{L^p([0,T],L^q)}\lesssim\|u_0\|_{\mathcal{H}^{\sigma}}+\|u_1\|_{\mathcal{H}^{\sigma-1}}.
$$
\end{theorem}

\medskip

In the case where the surface $M$ has a boundary, there is this slightly weaker result.

\medskip

\begin{theorem}
Let $(M,g)$ be a compact two-dimensional Riemannian manifold with boundary. Let $p\in(2,\infty]$ and $q\in[2,\infty)$ such that
$$
\frac{3}{p}+\frac{1}{q}\le\frac{1}{2}
$$
and consider $\sigma$ given by
$$
\frac{1}{p}+\frac{2}{q}=1-\sigma.
$$
Then the solution to
\begin{eqnarray*}
(\partial^2_t-\Delta_g)u & = & 0\\
(u,\partial_tu)|_{t=0} & = & (u_0,u_1)\in\mathcal{H}^{\sigma}\times\mathcal{H}^{\sigma-1}
\end{eqnarray*}
satisfies the bound
$$
\|u\|_{L^p([0,T],L^q)}\lesssim\|u_0\|_{\mathcal{H}^{\sigma}}+\|u_1\|_{\mathcal{H}^{\sigma-1}}.
$$
\end{theorem}

\bigskip

\subsection{Strichartz inequalities for wave equations with rough potentials}

While our proof of the Strichartz inequalities for the Schrödinger equation with white noise potential strongly relies on the deterministic result, this is not the case for the wave equation. In this case, we follow the approach from \cite{BLP} for which one has two mains ingredients, firstly a Weyl law for the Laplace-Beltrami operator and secondly $L^q$ bounds on its eigenfunctions. In particular, we treat at the same time the case with and without boundary here, the only difference being that one has weaker $L^q$ bounds on the eigenfunctions.

\medskip

An analogous Weyl law for the Anderson Hamiltonian was obtained in \cite{Mouzard}, see Proposition \ref{EingenBoundWeylLaw} in Section \ref{SubsectionAnderson}, and the analogue of the second part follows from the Strichartz inequalities for the Schrödinger group obtained in Section \ref{SectionStrichartzSchr}. Let $(e_n)_{n\ge1}$ be an orthonormal family of eigenfunctions of $H$ associated to $\big(\lambda_n(\Xi)\big)_{n\ge1}$. Since the eigenfuctions belong to the domain $\CD_\Xi$, they belong in particular to $L^\infty$ and we have the following bounds on its $L^q$-norm for $q\in(2,\infty)$. Recall that $(\lambda_n)_{n\ge1}$ are the eigenvalues of the Laplacian.

\medskip

\begin{proposition}
Let $q\in(2,\infty)$ and $M$ a compact surface without boundary. We have
$$
\|e_n\|_{L^q}\lesssim\sqrt{\lambda_n(\Xi)}^{\frac{1}{2}-\frac{1}{q}+\kappa}
$$
for any $\kappa>0$. In particular, this implies
$$
\|e_n\|_{L^q}\lesssim\big(1+\sqrt{\lambda_n}\big)^{\frac{1}{2}-\frac{1}{q}+\kappa} \lesssim(1+\sqrt{n})^{\frac{1}{2}-\frac{1}{q}+\kappa}.
$$
\end{proposition}

\medskip

\begin{proof}
We have
$$
\|e_n\|_{L^q}=\|e^{it\lambda_n}e_n\|_{L^p([0,1],L^q)}=\|e^{itH}e_n\|_{L^p([0,1],L^q)}
$$
with $(p,q)$ a Strichartz pair. For any $\kappa>0$, this gives
\begin{align*}
\|e_n\|_{L^q}&\lesssim\|e_n\|_{\CH^{\frac{1}{p}+\kappa}}\\
&\lesssim\|\sqrt{H}^{\frac{1}{p}+\kappa}e_n\|_{L^2}\\
&\lesssim\sqrt{\lambda_n(\Xi)}^{\frac{1}{p}+\kappa}
\end{align*}
using Proposition \ref{HSobolevBound} and
$$
\frac{1}{p}=\frac{1}{2}-\frac{1}{q}.
$$
Finally, Proposition \ref{EingenBoundWeylLaw} gives the bound
$$
\lambda_n(\Xi)\lesssim 1+\lambda_n
$$
and completes the proof.
\end{proof}

\medskip

Another important operator is the projection onto the eigenspaces of $H$. Let
$$\label{eqn:Pilambda}
\Pi_\lambda u:=\sum_{\lambda_n(\Xi)\in[\lambda,\lambda+1)}\langle u,e_n\rangle e_n
$$
for any $\lambda\ge0$. These spectral projectors satisfy the following bounds.

\medskip

\begin{proposition}\label{prop:eigen}
Let $\lambda\ge0$ and $q\in(2,\infty)$. We have
$$
\|\Pi_\lambda u\|_{L^q}\lesssim\sqrt{\lambda+1}^{\frac{1}{2}-\frac{1}{q}+\varepsilon}\|u\|_{L^2}
$$
for any $\varepsilon>0$.
\end{proposition}

\medskip

\begin{proof}
Consider $\lfloor H\rfloor$ the ``integer part'' of $H$ which is the self-adjoint operator defined by
$$
\lfloor H\rfloor e_n:=\lfloor\lambda_n(\Xi)\rfloor e_n
$$
for $n\ge1$. Then we have for any $\varepsilon>0$ the bound
$$
\|e^{it\lfloor H\rfloor}v\|_{L^p([0,1],L^q)}\lesssim\|v\|_{\CH^{\frac{1}{p}+\varepsilon}}
$$
which follows from the one for $H$, namely Theorem \ref{StrichartzSch}. Indeed, we have
$$
e^{it\lfloor H\rfloor}v-e^{itH}v=-i\int_0^te^{i(t-s)H}(H-\lfloor H\rfloor)e^{is\lfloor H\rfloor}v\drm s.
$$
and using Theorem \ref{StrichartzSch} and Corollary \ref{cor:inhom}, this gives
\begin{align*}
\|e^{it\lfloor H\rfloor}v\|_{L^p([0,1],L^q)}&\lesssim\|e^{itH}v\|_{L^p([0,1],L^q)}+\int_0^1\|e^{i(t-s)\lfloor H\rfloor}(H-\lfloor H\rfloor)e^{isH}\|_{L^p([0,1],L^q)}\drm s\\
&\lesssim\|v\|_{\CH^{\frac{1}{p}+\varepsilon}}+\int_0^1\|(H-\lfloor H\rfloor)e^{i(t-s)\lfloor H\rfloor}v\|_{\CH^{\frac{1}{p}+\varepsilon}}\drm s\\
&\lesssim\|v\|_{\CH^{\frac{1}{p}+\varepsilon}}
\end{align*}
for any $\varepsilon>0$ using that $\|H-\lfloor H\rfloor\|_{\CH^\beta\to\CH^\beta}$ is bounded by $1$ for $\beta<1$, which is true basically by construction together with Proposition \ref{HSobolevBound}, see also the proof of Proposition \ref{prop:floor}. Assuming that $\lambda\in\IN$, however the result follows directly in the same way by shifting $\lfloor H\rfloor$ for any  $\lambda\ge0$, we have
$$
\|e^{it\lfloor H\rfloor}\Pi_\lambda u\|_{L^p([0,1],L^q)}=\|e^{it\lambda}\Pi_\lambda u\|_{L^p([0,1],L^q)}=\|\Pi_\lambda u\|_{L^q}
$$
since the Weyl law guarantees that the number of eigenvalues in $[\lambda,\lambda+1)$ is finite. Thus we get using the Strichartz inequalities from Theorem \ref{StrichartzSch}
\begin{align*}
\|\Pi_\lambda u\|_{L^q}&\lesssim\|\Pi_\lambda u\|_{\CH^{\frac{1}{p}+\varepsilon}}\\
&\lesssim\sqrt{\lambda+1}^{\frac{1}{p}+\varepsilon}\|u\|_{L^2}\\
&\lesssim\sqrt{\lambda+1}^{\frac{1}{2}-\frac{1}{q}+\varepsilon}\|u\|_{L^2}
\end{align*}
using again Proposition \ref{HSobolevBound}.
\end{proof}

\medskip

As mentioned before, this is the point where there are slightly weaker results in the case of a surface with boundary. We use Theorem \ref{thm: strboundary} instead.

\medskip

\begin{proposition}\label{prop:eigenboundary}
Let $q\in(2,\infty)$ and $M$ a compact surface with boundary. We have
$$
\|e_n\|_{L^q}\lesssim\sqrt{\lambda_n(\Xi)}^{\frac{2}{3}-\frac{4}{3q}+\kappa}
$$
for any $\kappa>0$. In particular, this implies
$$
\|e_n\|_{L^q}\lesssim\big(1+\sqrt{\lambda_n}\big)^{\frac{2}{3}-\frac{4}{3q}+\kappa}\lesssim (1+\sqrt{n})^{\frac{2}{3}-\frac{4}{3q}+\kappa}.
$$
Moreover, for the operator $\Pi_\lambda$ we have
$$
\|\Pi_\lambda u\|_{L^q}\lesssim\sqrt{\lambda+1}^{\frac{2}{3}-\frac{4}{3q}+\kappa}\|u\|_{L^2}
$$
for any $\kappa>0$.
\end{proposition}

\medskip

Let $B$ be the operator defined by
$$
Be_n:=\lfloor\sqrt{\lambda_n(\Xi)}\rfloor e_n
$$
for any $n\ge1$. The following Proposition gives continuity estimates for the unitary groups associated to $\sqrt{H}$ and $B$ and bound the difference between the two operators.

\medskip

\begin{proposition}\label{prop:floor}
For any $\beta\in[0,1)$ and $t\in\IR$, we have
$$
\|e^{it\sqrt{H}}u\|_{\CH^\beta}\lesssim\|u\|_{\CH^\beta}
$$
and
$$
\|e^{itB}u\|_{\CH^\beta}\lesssim\|u\|_{\CH^\beta}.
$$
Moreover, the difference $B-\sqrt{H}$ is bounded on $\CH^\beta$ for any $\beta\in[0,1)$ and the difference between the groups is given by
$$
e^{itB}u-e^{it\sqrt{H}}=-i\int_0^te^{i(t-s)B}(\sqrt{H}-B)e^{is\sqrt{H}}\drm s.
$$
\end{proposition}

\medskip

\begin{proof}
We have
$$
\|e^{it\sqrt{H}}v\|_{L^2}\lesssim\|v\|_{L^2}.
$$
thus
$$
\|H^{\frac{\beta}{2}}e^{it\sqrt{H}}v\|_{L^2}=\|e^{it\sqrt{H}}H^{\frac{\beta}{2}}v\|_{L^2}\lesssim\|H^{\frac{\beta}{2}}v\|_{L^2}
$$
for any $\beta\in(0,\alpha)$. Using Proposition \ref{HSobolevBound}, this gives
$$
\|e^{it\sqrt{H}}v\|_{\CH^\beta}\lesssim\|v\|_{\CH^\beta}
$$
and the result for $e^{itB}$ follows from this. For the difference, $\|B-\sqrt{H}\|_{L^2\to L^2}$ is bounded by $1$ and we have
\begin{align*}
\|H^{\frac{\beta}{2}}(B-\sqrt{H})u\|_{L^2}&=\|(B-\sqrt{H})H^{\frac{\beta}{2}}u\|_{L^2}\\
&\le\|H^{\frac{\beta}{2}}u\|_{L^2}
\end{align*}
hence the boundedness of $B-\sqrt{H}$ on $\CH^\beta$. The result on the difference of the groups
$$
e^{itB}u-e^{it\sqrt{H}}=-i\int_0^te^{i(t-s)B}(\sqrt{H}-B)e^{is\sqrt{H}}\drm s
$$
follows with the same reasoning as in Lemma \ref{GroupDiff}.
\end{proof}

\medskip

We now have all the ingredients to prove of the Strichartz inequalities for the wave propagator associated to the Anderson Hamiltonian.

\medskip

\begin{theorem}\label{thm:strwaveM}
Let $M$ be a compact surface without boundary $(p,q)\in[2,\infty)^2$ and $0<\sigma<\alpha$ such that $p\le q$ and 
$$
\sigma=\frac{3}{2}-\frac{2}{p}+\frac{1}{q}.
$$
Then for any $\kappa>0$, we have the bound
$$
\big\|\cos(t\sqrt{H})u_0+\frac{\sin(t\sqrt{H})}{\sqrt{H}}u_1\big\|_{L^p([0,1],L^q)}\lesssim\|(u_0,u_1)\|_{\CH^{\sigma+\kappa}\times\CH^{\sigma-1+\kappa}}.
$$
\end{theorem}

\medskip

\begin{proof}
We start by proving the bound for $e^{itB}$ using the spectral decomposition
$$
e^{itB}u=\sum_{n\ge0}e^{itn}\Pi_nu
$$
and then bound the difference of the two groups. First, the condition $p\le q$ implies
$$
\|e^{itB}u\|_{L^p([0,1],L^q(M))}\le\|e^{itB}u\|_{L^q(M,L^p([0,1]))}
$$
hence it is enough to bound the right hand side. Using the Sobolev embedding in the time variable and the $L^q$ bound on the eigenvalues from Proposition \ref{prop:eigen}, we have
\begin{align*}
\|e^{itB}u\|_{L^q(M,L^p([0,1]))}^2&=\Big\|\|e^{itB}u\|_{L^p([0,1])}^2\Big\|_{L^{\frac{q}{2}}(M)}\\
&\lesssim\Big\|\|e^{itB}u\|_{\CH^{\frac{1}{2}-\frac{1}{p}}([0,1])}^2\Big\|_{L^{\frac{q}{2}}(M)}\\
&\lesssim\sum_{n\ge0}\Big\|\|e^{itn}\Pi_nu\|_{\CH^{\frac{1}{2}-\frac{1}{p}}([0,1])}^2\Big\|_{L^{\frac{q}{2}}(M)}\\
&\lesssim\sum_{n\ge0}\|e^{itn}\|_{\CH^{\frac{1}{2}-\frac{1}{p}}([0,1])}^2\|\Pi_nu\|_{L^q(M)}^2\\
&\lesssim\sum_{n\ge0}(n+1)^{1-\frac{2}{p}}(\sqrt{n}+1)^{1-\frac{2}{q}+2\kappa}\|\Pi_nu\|_{L^2}^2\\
&\lesssim\|\sqrt{H}^{\frac{3}{2}-\frac{2}{p}-\frac{1}{q}+\kappa}u\|_{L^2}^2\\
&\lesssim\|u\|_{\CH^{\frac{3}{2}-\frac{2}{p}-\frac{1}{q}+\kappa}}^2
\end{align*}
which gives the result for $B$. To obtain the proof for $\sqrt{H}$, we use
$$
e^{itB}u-e^{it\sqrt{H}}=-i\int_0^te^{i(t-s)B}(\sqrt{H}-B)e^{is\sqrt{H}}\drm s.
$$
Indeed, this gives
\begin{align*}
\|e^{it\sqrt{H}}u\|_{L^p([0,1],L^q)}&\lesssim\|e^{itB}u\|_{L^p([0,1],L^q)}+\int_0^1\|e^{i(t-s)B}(\sqrt{H}-B)e^{is\sqrt{H}}\|_{L^p([0,1],L^q)}\drm s\\
&\lesssim\|u\|_{\CH^{\sigma+\kappa}}+\int_0^1\|(\sqrt{H}-B)e^{i(t-s)B}u\|_{\CH^{\sigma+\kappa}}\drm s\\
&\lesssim\|u\|_{\CH^{\sigma+\kappa}}
\end{align*}
for any $\kappa>0$. The proof is directly completed from
$$
\cos(t\sqrt{H})=\frac{e^{it\sqrt{H}}+e^{-it\sqrt{H}}}{2}
$$
and
$$
\frac{\sin(\sqrt{H})}{\sqrt{H}}=\frac{e^{it\sqrt{H}}-e^{-it\sqrt{H}}}{2i\sqrt{H}}.
$$
\end{proof}

\medskip

Again, the inhomogeneous inequalities follow directly and we omit the proof.

\medskip

\begin{corollary}\label{cor:waveinhom}
Let $p,q,\sigma$ be as in Theorem \ref{thm:strwaveM}. Then we have the following bound
\begin{align*}
\big\|\int_{0}^{t}\frac{\sin\big((t-s)\sqrt{H}\big)}{\sqrt{H}}f(s)\big\|_{L^p([0,1],L^q)}\lesssim\int_0^1\|f(s)\|_{\CH^{\sigma-1+\kappa}}\drm s
\end{align*}
for $f\in L^1([0,1],\CH^{\sigma-1+\kappa}).$
\end{corollary}

\medskip

Moreover, we have the analogous result for surfaces with boundary which is proved analogously by using Proposition \ref{prop:eigenboundary} instead of Proposition \ref{prop:eigen}.

\medskip

\begin{theorem}\label{thm:strwavebound}
Let $M$ be a compact surface with boundary and $p,q\in[2,\infty)$ such that $p\le q$ and  
\begin{align*}
\sigma=\frac{5}{3}-\frac{2}{p}-\frac{4}{3q}\in(0,\alpha).
\end{align*}
Then for any $\kappa>0$, we have the bound
$$
\big\|\cos(t\sqrt{H})u_0+\frac{\sin(t\sqrt{H})}{\sqrt{H}}u_1+\int_{0}^{t}\frac{\sin((t-s)\sqrt{H})}{\sqrt{H}}v \big\|_{L^p([0,1],L^q)}\lesssim\|(u_0,u_1)\|_{\CH^{\sigma+\kappa}\times\CH^{\sigma-1+\kappa}}+\|v\|_{L^1({[0,1]},\CH^{\sigma-1+\kappa})}
$$
for initial data $(u_0,u_1)\in\mathcal{H}^\sigma\times\mathcal{H}^{\sigma-1}$ and inhomogeneity $v\in L^1({[0,1]},\CH^{\sigma-1+\kappa})$.
\end{theorem}

\bigskip

\subsection{Local well-posedness for the multiplicative cubic stochastic wave equation}

Now we use the results from the previous Section to prove local well-posedness of stochastic multiplicative wave equations of the form
$$\label{eqn:HWave}
\left|\begin{array}{ccc}
\partial^2_tu+Hu&=&-u|u|^2\\
(u,\partial_tu)|_{t=0}&=&(u_0,u_1)
\end{array}\right.
$$
in a low-regularity regime on general two-dimensional surfaces with or without boundary. While we have the classical Sobolev embedding
$$
\CH^\nu\hookrightarrow L^{\frac{2}{1-\nu}}
$$
for $\nu\in[0,1)$, we also make use of the following dual Sobolev bound
$$\label{eqn:dualSob}
\forall\sigma\in(0,1],\quad L^{\frac{2}{2-\sigma}}\hookrightarrow\CH^{\sigma-1}
$$ 
which is true on general manifolds, see for example the book by Aubin \cite{Aubin}. Using this, we make a preliminary computation meant to show how far we get by using \textit{only the Sobolev embedding result}. Then we will see how the bounds in Theorem \ref{thm:strwaveM} give better results on general manifolds. We first rewrite the equation under the mild formulation
$$
u(t)=\cos(t\sqrt{H})u_0+\frac{\sin(t\sqrt{H})}{\sqrt{H}}u_1+\int_{0}^{t}\frac{\sin\big((t-s)\sqrt{H}\big)}{\sqrt{H}}u(s)^3\drm s.
$$
Then apply the dual Sobolev bound for $\sigma\in(0,1]$ and $p=\frac{2}{2-\sigma}\in(1,2]$ to get
\begin{align*}
\|u(t)\|_{\mathcal{H}^{\sigma}}\lesssim&\|u_0\|_{\mathcal{H}^{\sigma}}+\|u_1 \|_{\mathcal{H}^{\sigma-1}}+\|u^3\|_{L^1([0,t],\mathcal{H}^{\sigma-1})}\\
\lesssim&\|u_0\|_{\mathcal{H}^{\sigma}}+\|u_1\|_{\mathcal{H}^{\sigma-1}} +\|u^3\|_{L^1([0,t],L^{p})}\\
\lesssim& \|u_0\|_{\mathcal{H}^{\sigma}}+\|u_1\|_{\mathcal{H}^{\sigma-1}}+\|u\|_{L^\infty([0,t],L^{\frac{2}{1-\sigma}})}\|u\|^2_{L^2([0,t],L^{4})},
\end{align*}
having applied Hölder with $\frac{1}{2}+\frac{1-\sigma}{2}=\frac{2-\sigma}{2}.$ Finally, the Sobolev embedding gives
\begin{align*}
\|u(t)\|_{\CH^\sigma}\lesssim&\|u_0\|_{\mathcal{H}^{\sigma}}+\|u_1\|_{\mathcal{H}^{\sigma-1}}+\|u\|_{L^\infty([0,t],\mathcal{H}^{\sigma})}\|u\|^2_{L^2([0,t],\mathcal{H}^{\frac{1}{2}})}.
\end{align*}
This can then lead to a solution by fixed point by choosing $\sigma\ge\frac{1}{2}.$ Clearly this is can be improved by using more subtle bounds than the Sobolev embedding. The Strichartz inequalities from the previous section allow us to get local well-posedness below, this is the content of the following Theorems; As before we separately state the cases of surfaces without boundary, with boundary which are proved in precisely the same way, just using Theorems \ref{thm:strwaveM} and \ref{thm:strwavebound} respectively.

\medskip

\begin{theorem}
Let $M$ be a compact surface without boundary and $\sigma\in(\frac{1}{4},\frac{1}{2})$ and $\delta>0$ sufficiently small. Then for any initial data $(u_0,u_1)\in \mathcal{H}^\sigma\times\mathcal{H}^{\sigma-1}$ there exists a time $T>0$ depending on the data such that there exists a unique solution
$$
u\in C\big([0,T],\mathcal{H}^\sigma\big)\cap L^{\frac{2}{1-\delta}}\big([0,T],L^{4}\big)
$$
to the mild formulation of the multiplicative cubic stochastic wave equation. Moreover, the solution depends continuously on the initial data $(u_0,u_1)$.
\end{theorem}

\medskip

\begin{proof}
As usual, this is proved in a standard way using the Banach fixed point Theorem. Define the map
$$
\Psi (u)(t):=\cos(t\sqrt{H})u_0+\frac{\sin(t\sqrt{H})}{\sqrt{H}}u_1+\int_{0}^{t}\frac{\sin\big((t-s)\sqrt{H}\big)}{\sqrt{H}}u(s)^3\drm s.
$$		
For $t>0$, we have as above
\begin{align*}
\|u(t)\|_{\mathcal{H}^{\sigma}}&\lesssim \|u_0\|_{\mathcal{H}^{\sigma}}+\|u_1\|_{\mathcal{H}^{\sigma-1}}+\|u\|_{L^\infty([0,t],\mathcal{H}^{\sigma})}\|u\|^2_{L^2([0,t],L^{4})}\\
&\lesssim \|u_0\|_{\mathcal{H}^{\sigma}}+\|u_1\|_{\mathcal{H}^{\sigma-1}}+t^\delta\|u\|_{L^\infty([0,t],\mathcal{H}^{\sigma})}\|u\|^2_{L^\frac{2}{1-\delta}([0,t],L^{4})}
\end{align*}
using Hölder inequality in the last line for $\delta\in(0,1)$. We then apply Theorem \ref{thm:strwaveM} with $p=\frac{2}{1-\delta}$ and $q=4$ and obtain
\begin{align*}
\big\|\Psi(u)\big\|_{L^{\frac{2}{1-\delta}}\big([0,T],L^{4}\big)}&\lesssim\|u_0\|_{\CH^{\frac{3}{2}-(1-\delta)- \frac{1}{4}+\kappa}}+\|u_1\|_{\CH^{\frac{1}{2}-(1-\delta)-\frac{1}{4}+\kappa}}+\|u^3\|_{L^1([0,T],\CH^{\frac{1}{2}-(1-\delta)-\frac{1}{4}+\kappa})}\\
&\lesssim\|u_0\|_{\CH^\sigma}+\|u_1\|_{\CH^{\sigma-1}}+\|u^3\|_{L^1([0,T],\CH^{\sigma-1})}
\end{align*}
using that $\sigma>\frac{1}{4}$ and $\delta<\sigma-\frac{1}{4}$ gives $\frac{3}{2}-(1-\delta)-\frac{1}{4}+\kappa\le\sigma$ for $\kappa>0$ small enough. Finally, we get
\begin{align*}
\big\|\Psi(u)\big\|_{L^{\frac{2}{1-\delta}}\big([0,T],L^{4}\big)}
&\lesssim\|u_0\|_{\mathcal{H}^{\sigma}}+\|u_1\|_{\mathcal{H}^{\sigma-1}}+T^\delta\|u\|_{L^\infty([0,T],\mathcal{H}^{\sigma})}\|u\|^2_{L^\frac{2}{1-\delta}([0,T],L^{4})}
\end{align*}
as above. Thus we can get a fixed point in
$$
C\big([0,T],\mathcal{H}^\sigma\big)\cap L^{\frac{2}{1-\delta}}\big([0,T],L^{4}\big)
$$
in the usual way for $T>0$ small enough.
\end{proof}

\medskip

In a completely analogous way we get the following result for the case of surfaces with boundary using the Strichartz estimates from Theorem \ref{thm:strwavebound}.

\medskip

\begin{theorem}
Let $M$ be a compact surface with boundary and $\sigma\in(\frac{1}{3},\frac{1}{2})$ and $\delta>0$ sufficiently small. Then for any initial data $(u_0,u_1)\in \mathcal{H}^\sigma\times\mathcal{H}^{\sigma-1}$ there exists a time $T>0$ depending on the data such that there exists a unique solution
$$
u\in C\big([0,T],\mathcal{H}^\sigma\big)\cap L^{\frac{2}{1-\delta}}\big([0,T],L^{4}\big)
$$
to the mild formulation of the multiplicative cubic stochastic wave equation. Moreover, the solution depends continuously on the initial data $(u_0,u_1)$.
\end{theorem}

\bigskip





\bigskip

\noindent \textcolor{gray}{$\bullet$} A. Mouzard --  Univ. Rennes, CNRS, IRMAR - UMR 6625, F-35000 Rennes, France   \\
{\it E-mail}: antoine.mouzard@ens-rennes.fr

\noindent \textcolor{gray}{$\bullet$} I. Zachhuber --  , FU Berlin, D-14195 Berlin, Germany\\
{\it E-mail}: immanuel.zachhuber@fu-berlin.de

\end{document}